\documentclass[11pt]{article}
\usepackage{amsmath,amssymb,amsthm}
\usepackage{graphics}

\usepackage{latexsym}		
\usepackage{graphicx}		
\usepackage{color}
\usepackage{subfigure}
\usepackage{amsfonts}
\usepackage{mathtools}

\newcommand{\beq}{\begin{equation}}
\newcommand{\eeq}{\end{equation}}
\newcommand{\bdm}{\begin{displaymath}}
\newcommand{\edm}{\end{displaymath}}
\newcommand{\beqa}{\begin{eqnarray}}
\newcommand{\eeqa}{\end{eqnarray}}
\newcommand{\beqas}{\begin{eqnarray*}}
\newcommand{\eeqas}{\end{eqnarray*}}

\newcommand{\BU}{{\bf U}}

\newcommand{\bb}{{\bf b}}

\newcommand{\bsigma}{\boldsymbol\sigma}

\newcommand{\bPsi}{\boldsymbol\Psi}
\newcommand{\bnu}{\boldsymbol\nu}

\newcommand{\brho}{\boldsymbol\rho}
\newcommand{\betab}{\boldsymbol\eta}
\newcommand{\balpha}{\boldsymbol\alpha}
\newcommand{\bbeta}{\boldsymbol\beta}
\newcommand{\bgamma}{\boldsymbol\gamma}
\newcommand{\bkappa}{\boldsymbol\kappa}

\newcommand{\bxi}{\boldsymbol\xi}

\newcommand{\br}{{\bf r}}
\newcommand{\bu}{{\bf u}}
\newcommand{\bv}{{\bf v}}
\newcommand{\bw}{{\bf w}}
\newcommand{\bx}{{\bf x}}
\newcommand{\by}{{\bf y}}
\newcommand{\bz}{{\bf z}}

\newcommand{\bbR}{\mathbb{R}}

\newcommand{\bRn}{{\mathbb{R}^n}}

\newcommand{\CC}{{\cal C}}
\newcommand{\CD}{{\cal D}}

\newcommand{\CG}{{\cal G}}

\newcommand{\CL}{{\cal L}}

\newcommand{\TTT}{{\bf T}}
\newcommand{\ttt}{{\sigma}}
\newcommand{\CCC}{{\bf C}}
\newcommand{\SSS}{{\bf S}}
\newcommand{\KKK}{{\bf K}}

\newcommand{\bzeta}{\boldsymbol\zeta}

\newcommand{\Us}{{[C_c^{\infty}(\bRn\times\bRn)]^k}}
\newcommand{\Vs}{{C_c^{\infty}(\bRn)}}
\newcommand{\Usp}{{[D'(\bRn\times\bRn)]^k}}
\newcommand{\Vsp}{{D'(\bRn)}}
\newcommand{\UVs}{{\Us\times \Vs}}
\newcommand{\UVsp}{{\big(\Us\times \Vs\big)'}}

\newtheorem{defn}{Definition}
\newtheorem{thm}{Theorem}
\newtheorem{cor}{Corollary}
\newtheorem{lem}{Lemma}

\newtheorem{prop}{Proposition}

\numberwithin{equation}{section}

\begin{document}
\nocite{*}
\title{A Generalized Nonlocal Calculus with Application to the Peridynamics Model for Solid Mechanics}
\author{Bacim Alali, Kuo Liu, and Max Gunzburger\\
\
\footnotesize{Department of Scientific Computing, Florida State University, Tallahassee, FL}
\thanks{Research and preparation of the paper was partially supported by the US Department of Energy grant number DE-SC0004970
 and by the US National Science Foundation grant number DMS-1013845.}}
\date{}
\maketitle
%
%
%
%
%
%
\begin{abstract}
A nonlocal vector calculus was introduced in \cite{dglz-nlc} that has proved useful for the  analysis of 
the peridynamics model of nonlocal mechanics  and nonlocal diffusion models. A generalization is developed that provides a more general setting for the nonlocal vector calculus that is independent of particular nonlocal models.
It is shown that general nonlocal calculus operators are integral operators with specific integral kernels. 
General nonlocal calculus properties are developed, including nonlocal integration by parts formula and Green's identities. 
The nonlocal vector calculus introduced in \cite{dglz-nlc} is shown to be recoverable from the general formulation as a special example. 
This special nonlocal vector calculus is used to reformulate the peridynamics equation of motion in terms of the nonlocal 
gradient operator and its adjoint.
A new example of nonlocal vector calculus operators is introduced, which shows the potential use of the 
general formulation for general nonlocal models. 
\end{abstract}
%
%
%
%

\section{Introduction}
In recent years, nonlocal continuum models have been developed for several large-scale phenomena. 
Examples include the peridynamics formulation
for solid mechanics \cite{Silling2000,Silling2007} and nonlocal diffusion \cite{dglz-nld}.
These nonlocal continuum models are described through integral equations  in contrast to their classical 
local continuum counterparts which are given by partial differential equations. 
A key connection between the peridynamics model and classical elasticity and between the nonlocal diffusion model 
and classical diffusion is that these nonlocal models have been shown to converge, under certain conditions, 
 to their local counterparts in the limit of vanishing nonlocality \cite{emmrich2007well,dglz-ps,dglz-nld}.
Another connection between these local and nonlocal models is given through a {\it nonlocal vector calculus} that is 
introduced and developed in \cite{dglz-nlc}. The nonlocal vector calculus introduces integral operators that mimic the roles of
the divergence, gradient, and other vector calculus operators.
Specifically, the nonlocal divergence of a vector-valued function $\bnu(\bx,\by)$ is defined as \cite{dglz-nlc}
\begin{equation}
\label{sdiv0}
 (\CD_{\balpha}\bnu)(\bx) = \int\big(\bnu(\bx,\by) + \bnu(\by,\bx)\big)\cdot\balpha(\bx,\by)\,d\by,
\end{equation}
where the kernel $\balpha$ is an antisymmetric vector-valued  function, i.e, $\balpha(\bx,\by)=-\balpha(\by,\bx)$.
In addition, the action of the adjoint operator $\CD^\ast_{\balpha} $ on a scalar function $u(\bx)$ is given by 
\begin{equation}
\label{sdivadj0}
 (\CD_{\balpha}^\ast u)(\bx,\by) = -\big(u(\by) - u(\bx)\big)\balpha(\bx,\by).
\end{equation}
Moreover, for a scalar function $\eta(\bx,\by)$ and a vector-valued function $\bu(\bx)$, 
the nonlocal gradient operator $\CG_{\balpha}$ and its adjoint $\CG_{\balpha}^{\ast}$ are defined by
\begin{eqnarray}
\label{sgrad0}
 (\CG_{\balpha}\eta)(\bx) &=& \int\big(\eta(\by,\bx) + \eta(\bx,\by)\big)\balpha(\bx,\by)\, d\by,\\
\label{sgradadj0}
 (\CG_{\balpha}^* \bu)(\bx,\by) & =& -\big(\bu(\by)-\bu(\bx)\big)\cdot\balpha(\bx,\by).
\end{eqnarray}
Using these nonlocal operators, with specific choices of $\balpha$, it is shown in \cite{dglz-nld} that 
\[
 \dot{u}+\CD_{\balpha}(\CD_{\balpha}^\ast u)=b
\]
is a nonlocal diffusion equation. In addition, the linear peridynamics equation \cite{Silling2007}
 \begin{equation}
\label{peridyn_0}
    \ddot{\bu}=\mathcal{L}\bu+\bb,
 \end{equation}
where $\mathcal{L}$ is given by \eqref{lpd3}, 
can be written, using nonlocal vector calculus operators \cite{dglz-ps}, as
\begin{equation}
\label{peridyn_L_old}
 \mathcal{L}\bu=-\CD_{\balpha}(c_{1}' \omega \;(\CD_{\balpha}^{\ast}\bu)^T)-
\CD_{\balpha}^{\omega}(c_{2}' \mbox{tr} ({\CD_{\balpha}^{\omega}}^{\ast}\bu) I),
\end{equation}
where $I$ is the identity matrix, $c_1 ', c_2 '$ are material properties, $\omega$ a weight function, 
and  $\CD_{\balpha}^{\omega}$, ${\CD_{\balpha}^{\omega}}^{\ast}$
are weighted versions of $\CD_{\balpha}$, $\CD_{\balpha}^{\ast}$, respectively; see \cite{dglz-ps} for details.

In this work, we show that the linear peridynamics operator $\mathcal{L}$ has a simpler expression in terms of nonlocal vector calculus
operators. In Theorem \ref{thmpdandno}, we show that, for an appropriate choice of the integral kernel $\balpha$, the peridynamics
operator $\mathcal{L}$ in \eqref{peridyn_0} can be cast as
\begin{equation}
\label{peridyn_L_new}
 \mathcal{L}\bu=-\CG_{\balpha}(c_{1} \CG_{\balpha}^{\ast}\bu)-\CG_{\balpha}(c_{2} \overline{\CG_{\balpha}^{\ast}}\bu),
\end{equation}
where $c_1, c_2$ are scalars, and $\overline{\CG_{\balpha}^{\ast}}$ is an average of $\CG_{\balpha}^{\ast}$ defined by
\[
 (\overline{\CG_{\balpha}^*} \bu)(\bx)  = -\int \big(\bu(\by)-\bu(\bx)\big)\cdot\balpha(\bx,\by)d\by.
\]
This  new expression  for $\mathcal{L}$ given in \eqref{peridyn_L_new} bears a closer resemblance to the Navier operator of linear elasticity.

Given the fact that the nonlocal  calculus operators given by \eqref{sdiv0}--\eqref{sgradadj0} 
 mimic the  differential calculus operators
in the setting of nonlocal diffusion and peridynamics models, 
one may ask whether these operators are the only nonlocal integral operators that do so.
In this work, we provide a general mathematical setting for the existence of nonlocal integral operators that resemble the differential 
calculus operators independent of particular nonlocal models. In Section~\ref{gennlc}, we show that a nonlocal  operator that 
resembles\footnote{The resemblance of nonlocal divergence to local divergence is made precise in Section~\ref{gennlc}.} the divergence operator, for instance, must be of the general form
\begin{equation}
\label{general_D}
(\CD \bnu )(\bx) = 
\int\int\bkappa(\bx,\by,\bz)\cdot\bnu(\by,\bz)\,d\bz d\by,
\end{equation}
for some kernel $\bkappa$ that satisfies
\begin{equation}
\label{div_cond}
\int\bkappa(\bx,\by,\bz)\,d\bx = 0
\qquad\mbox{for a.e. $\by,\bz$.}
\end{equation}
We refer to the operator $\CD$ in \eqref{general_D}--\eqref{div_cond} as general nonlocal divergence. We introduce general nonlocal operators including a nonlocal gradient, nonlocal curl,  and nonlocal Laplacian. General nonlocal calculus theorems and identities such as  nonlocal integration by parts formulas and Green's identities are developed.

We show in Section~\ref{sec_special} that the nonlocal divergence $\CD_{\balpha}$ in \eqref{sdiv0} can be recovered from \eqref{general_D} for a specialized kernel $\bkappa=\bkappa(\balpha)$. The other nonlocal operators in \eqref{sdivadj0}--\eqref{sgradadj0} are also shown to follow from the general formulation of the nonlocal calculus.

In Section~\ref{sec_conclusion}, we provide a new example for  nonlocal  calculus operators. Specifically, we show that the operator defined by
\begin{equation}
\label{sdiv_beta}
 (\CD_{\bbeta}\bnu)(\bx) = \int\big(\bnu(\by,\bx) - \bnu(\bx,\by)\big)\cdot\bbeta(\bx,\by)\,d\by,
\end{equation}
where the  kernel $\bbeta$ is a symmetric vector-valued  function, is a nonlocal divergence operator. The operator $\CD_{\bbeta}$ is a special case of  \eqref{general_D} for a specific kernel $\bkappa=\bkappa(\bbeta)$.
It is anticipated that nonlocal calculus operators, such as $\CD_{\bbeta}$ in \eqref{sdiv_beta}, will be useful for the analysis of new  nonlocal models.

This article is organized as follows. Section~\ref{gennlc} introduces the general formulation for the nonlocal vector calculus. General nonlocal calculus theorems, identities, and regularity results for nonlocal operators are derived. Section~\ref{sec_special} focuses on the special case of nonlocal calculus operators defined in \eqref{sdiv0}--\eqref{sgradadj0}. An application to the peridynamics model of solid mechanics is discussed in Section~\ref{sec_peridynamics}. Conclusion remarks and discussion of a new example of nonlocal calculus operators are provided in Section~\ref{sec_conclusion}.

\section{A generalized nonlocal calculus}\label{gennlc}

For the spaces $\Us$ and $\Vs$ and the corresponding dual spaces \\ $\Usp$ and $\Vsp$ with $k=1,2$ or $3$, we have the duality parings 
\[
\begin{aligned}
<\bnu,\bgamma>_{\Us,\Usp} &= \int_{\bRn}\int_{\bRn}\bnu(\bz,\by)\cdot\bgamma(\bz,\by)\,d\bz d\by,
\\&
\hspace*{-1cm}\forall \,\bnu\in \Us,\,\,\bgamma\in \Usp
\end{aligned}
\]
and 
$$
<v,u>_{\Vs,\Vsp} = \int_{\bRn}u(\bx)v(\bx)\,d\bx
\qquad\forall\, v\in \Vs,\,\,u\in \Vsp.
$$
For the product space $\UVs$ and its dual\\ $\UVsp$, the duality paring is given by
$$
\begin{aligned}
&<\bkappa,\,\bsigma>_{\UVs, \UVsp} 
\\&\qquad\qquad= \int_{\bRn}\int_{\bRn}\int_{\bRn}\bsigma(\bx,\by,\bz)\cdot\bkappa(\bx,\by,\bz)\,d\bz d\by d\bx,
\\&
\forall \,
\bsigma\in \UVsp,\,\, \bkappa\in\UVsp.
\end{aligned}
$$

Let $\CD: \Us\to \Vsp$ denote a linear and continuous operator. Then, by the Schwartz Kernel Theorem, there exists a unique $\bkappa\in$\\ $\UVsp$ such that
$$
\begin{aligned}
&<v,\CD \bnu>_{\Vs,\Vsp}\\
& \qquad\qquad= <v\otimes\bnu,\bkappa>_{\UVs,\UVsp},\\
&\forall\,v\in \Vs,\,\,\bnu\in \Us
\end{aligned}
$$
or, using the definitions of the duality pairings,
$$
\begin{aligned}
\int_{\bRn}v(\bx)(\CD \bnu )(\bx)\,d\bx
=\int_{\bRn}\int_{\bRn}\int_{\bRn} &v(\bx)\bnu(\by,\bz)\cdot\bkappa(\bx,\by,\bz)\,d\bz d\by d\bx
\\&
\forall\,v\in \Vs,\,\,\bnu\in \Us.
\end{aligned}
$$
The arbitrariness of $v\in \Vs$ implies that $\CD \bnu \in \Vsp$ is given by
\begin{equation}\label{sktd}
(\CD \bnu )(\bx) = 
\int_{\bRn}\int_{\bRn}\bkappa(\bx,\by,\bz)\cdot\bnu(\by,\bz)\,d\bz d\by
\qquad\mbox{for almost all $\bx\in \bRn$.}
\end{equation}

We seek an operator $\CD$ that satisfies a divergence-like theorem which we now describe. For $\bnu\in [L^1(\bRn\times\bRn)]^k$, let $\psi_{\bnu}\in[L^1(\bRn\times\bRn)]^k$ be such that
\begin{subequations}\label{psiass}
\begin{align}
&\mbox{$\psi_{\bnu}$ is linear in $\bnu$}\label{psiassa}\\
&\mbox{$\psi_{\bnu}$ is antisymmetric, i.e., $\psi_{\bnu}(\bx,\by)=-\psi_{\bnu}(\by,\bx)$ for all $\bx,\by\in\bRn$}.\label{psiassb}
\end{align}
\end{subequations}
For any $\bx\in\bRn$ and $\widetilde\Omega\subset\bRn$, $\int_{\widetilde\Omega}\psi_{\bnu}(\bx,\by)\,d\by$ represents the nonlocal flux density at $\bx$ into $\widetilde\Omega$; see \cite{dglz-nlc} for details. The operator $\CD$ and the flux density $\int_{\widetilde\Omega}\psi_{\bnu}(\bx,\by)\,d\by$ are required to satisfy the nonlocal ``divergence'' theorem\footnote{In words, \eqref{divthm} states that the integral of the nonlocal divergence of $\bnu$ over any domain $\Omega\subset\bRn$ is equal to the flux of $\bnu$ exiting from $\Omega$ into the complement domain $\bRn\setminus\Omega$. This is made clear by noting that, due to the antisymmetry of $\psi_{\bnu}(\bx,\by)$, \eqref{divthm} can be rewritten as
$$
\int_\Omega (\CD\bnu)(\bx)\,d\bx = 
\int_\Omega \int_{\bRn\setminus\Omega}\psi_{\bnu}(\bx,\by)\,d\by d\bx
\qquad
\forall\,\bnu\in \Us,\,\, \Omega\subset\bRn.
$$}
\begin{equation}\label{divthm}
\int_\Omega (\CD\bnu)(\bx)\,d\bx = 
\int_\Omega \int_\bRn\psi_{\bnu}(\bx,\by)\,d\by d\bx
\qquad
\forall\,\bnu\in \Us,\,\, \Omega\subset\bRn.
\end{equation}

From \eqref{divthm} and the arbitrariness of $\Omega$, we obtain
\begin{equation}\label{divthm2}
(\CD\bnu)(\bx) = \int_\bRn\psi_{\bnu}(\bx,\by)\,d\by
\qquad\mbox{for a.e. $\bx\in\bRn$.}
\end{equation}
From \eqref{sktd} and \eqref{divthm2}, we obtain
\begin{equation}\label{divthm3}
\int_\bRn\psi_{\bnu}(\bx,\by)\,d\by
=\int_{\bRn}\int_{\bRn}\bkappa(\bx,\by,\bz)\cdot\bnu(\by,\bz)\,d\bz d\by \qquad\mbox{for a.e. $\bx\in\bRn$.}
\end{equation}
Note that here $\bkappa$ is fixed whereas $\bnu$ is arbitrary.

\begin{lem}\label{lemrho}
The kernel $\bkappa(\bx,\by,\bz)$ satisfies
\begin{equation}\label{divthm4}
\int_{\bRn}\bkappa(\bx,\by,\bz)\,d\bx = 0
\qquad\mbox{for a.e. $\by,\bz\in\bRn$.}
\end{equation}
\end{lem}

{\em Proof.}
From \eqref{divthm} and the antisymmetry of $\psi_{\bnu}(\bx,\by)$, we have
\begin{equation}\label{divthm5}
\int_\bRn (\CD\bnu)(\bx)\,d\bx = 
\int_\bRn \int_\bRn\psi_{\bnu}(\bx,\by)\,d\by d\bx
=-\int_\bRn \int_\bRn\psi_{\bnu}(\by,\bx)\,d\bx d\by
=0.
\end{equation}
Then, from \eqref{divthm2}, \eqref{divthm3}, and \eqref{divthm5}, we have
\begin{eqnarray*}
\int_{\bRn}\int_{\bRn}\int_{\bRn}\bkappa(\bx,\by,\bz)\cdot\bnu(\by,\bz)\,d\bz d\by d\bx = \int_\bRn (D\bnu)(\bx)\,d\bx =0,\\
\quad\forall\,\bnu\in \Us .
\end{eqnarray*}
Therefore,
$$
\int_{\bRn}\int_{\bRn}\Big(\int_{\bRn}\bkappa(\bx,\by,\bz)\,d\bx\Big)\cdot\bnu(\by,\bz)\,d\bz d\by =0,
\qquad \forall\,\bnu\in \Us 
$$
which implies \eqref{divthm4}.
\hfill$\Box$

We refer to a kernel $\bkappa(\bx,\by,\bz)$ satisfying \eqref{divthm4} as a {\em divergence kernel.}

From \eqref{divthm2} and \eqref{divthm3}, we are led to the following definition of a nonlocal divergence operator.

\begin{defn}[{\bf Nonlocal divergence operator}]\label{ddiv}
The action of the {\em nonlocal divergence operator} $\CD: \Us\to \Vsp$ on any vector-valued function $\bnu\in \Us$ is given by\footnote{A similar definition to \eqref{divop} was given in \cite{dglz-nlc}. However, there, the central requirement \eqref{divthm4} was not discussed nor was the development of the full nonlocal vector calculus associate with \eqref{divop}.}
\begin{equation}\label{divop}
(\CD\bnu)(\bx) = 
\int_{\bRn}\int_{\bRn}\bkappa(\bx,\by,\bz)\cdot\bnu(\by,\bz)\,d\bz d\by \qquad\mbox{for a.e. $\bx\in \bRn$,}
\end{equation}
where $\bkappa\in\UVsp$ satisfies \eqref{divthm4}.
\hfill$\Box$
\end{defn}

The adjoint operator $\CD^\ast$ corresponding to the nonlocal divergence operator $\CD$ is defined through the relation
\begin{equation}\label{adjdef}
\begin{aligned}
  < u,\CD\bnu >_{\Vs,\Vsp} = <\bnu, &\CD^\ast u>_{\Us,\Usp}
  \\&\forall\, \bnu\in \Us,\,\,u\in \Vs.
\end{aligned}
\end{equation}

\begin{prop}[\bf Adjoint operator]\label{propadj}
Corresponding to the nonlocal divergence operator $\CD: \Us\to \Vsp$, we have the {\em adjoint operator} $\CD^\ast: \Vs\to \Usp$ whose action on any scalar-valued function $u\in \Vs$ is given by\footnote{With $\CD^\ast$ being the adjoint of the nonlocal divergence operator $\CD$, one can identify $-\CD^\ast$ as a nonlocal gradient operator.}
\begin{equation}\label{adjop}
   (\CD^\ast u)(\bx,\by) = \int_\bRn u(\bz) \bkappa(\bz,\bx,\by)\,d\bz \qquad\mbox{for a.e. $\bx,\by\in \bRn$}.
\end{equation}
\end{prop}

{\em Proof.}
From \eqref{divop} and \eqref{adjdef} we have
\begin{equation}\label{adjpf}
\begin{aligned}
\int_{\bRn}\Big(\int_{\bRn}\int_{\bRn}&\bkappa(\bx,\by,\bz)\cdot\bnu(\by,\bz)\,d\bz d\by\Big) u(\bx)\,d\bx
\\&= \int_{\bRn}\int_{\bRn}
(\CD^\ast u)(\bx,\by)\cdot\bnu(\bx,\by)\,d\by d\bx.
\end{aligned}
\end{equation}
After switching the dummy variables $\bx$ and $\bz$ and then $\bx$ and $\by$ in the left-hand side, we have
$$
\begin{aligned}
\int_{\bRn}\int_{\bRn}\int_{\bRn}&\bkappa(\bx,\by,\bz)\cdot\bnu(\by,\bz)u(\bx)\,d\bz d\by d\bx
\\&=
\int_{\bRn}\int_{\bRn}\int_{\bRn}\bkappa(\bz,\by,\bx)\cdot\bnu(\by,\bx)u(\bz)\,d\bx d\by d\bz
\\&=
\int_{\bRn}\int_{\bRn}\int_{\bRn}\bkappa(\bz,\bx,\by)\cdot\bnu(\bx,\by)u(\bz)\,d\by d\bx d\bz
\\&=
\int_{\bRn}\int_{\bRn}\Big(\int_{\bRn}\bkappa(\bz,\bx,\by)u(\bz)\,d\bz\Big)\cdot\bnu(\bx,\by)\,d\by d\bx. 
\end{aligned}
$$
Then, because $\bnu\in \Us$ is arbitrary, we obtain \eqref{adjop} from \eqref{adjpf}.
\hfill$\Box$


%
%
%
%
\subsection{Regularity of $\CD$ and $\CD^\ast$}

In Definition \ref{ddiv} and Proposition \ref{propadj}, we assume that $\bnu(\bx,\by)\in\Us$, $u\in \Vs$, and $\bkappa\in\UVsp$. In fact, the nonlocal divergence operator $\CD$ and its adjoint operator $\CD^*$ can be defined for functions having much less smoothness, as the next proposition shows.

\begin{prop}[\bf Regularity of $\CD$ and $\CD^*$]\label{propreg}
Let $1\leq p \leq \infty$ and $1\leq q \leq \infty$ with $\frac{1}{p} + \frac{1}{q} = 1$ and assume that $\bkappa\in [L^p(\bRn\times\bRn\times\bRn)]^k$. Then,
\begin{subequations}
\begin{align} 
 \CD&:\;\bnu\in [L^q(\bRn\times\bRn)]^k\longmapsto \CD\bnu \in L^p(\bRn) \label{drega}
\\
\CD^*&:\; u\in L^q(\bRn)\longmapsto \CD^*u \in [L^p(\bRn\times\bRn)]^k.\label{dregb}
\end{align}
\end{subequations}
In particular, $[L^2(\bRn\times\bRn)]^k \xmapsto{\CD} L^2(\bRn) \xmapsto{\CD^*} L^2[(\bRn\times\bRn)]^k$. Moreover, $\CD$ and $\CD^*$ are bounded operators on $[L^2(\bRn\times\bRn)]^k$ and $L^2(\bRn)$, respectively. 
\end{prop}

{\em Proof.}
Letting $\bnu\in [L^q(\bRn\times\bRn)]^k$, by Minkowski's integral inequality, we have 
$$
\begin{aligned}
 \|\CD\bnu\|_{L^p(\bRn)} &= \Big(\int_{\bRn} \Big|\int_{\bRn}\int_{\bRn}\bnu(\by,\bz)\cdot\bkappa(\bx,\by,\bz)\,d\bz d\by\Big|^p d\bx \Big)^{1/p} \\
                         &\leq \int_{\bRn}\int_{\bRn} \Big( \int_{\bRn}|\bnu(\by,\bz)\cdot\bkappa(\bx,\by,\bz)|^p d\bx\Big)^{1/p}\,d\bz d\by \\
                         &\leq \int_{\bRn}\int_{\bRn} |\bnu(\by,\bz)| \Big( \int_{\bRn}|\bkappa(\bx,\by,\bz)|^p\,d\bx\Big)^{1/p} d\bz d\by.
\end{aligned}
$$
Applying the H\"{o}lder inequality to the last inequality, we have 
\begin{equation}\label{dregc}
 \|\CD\bnu\|_{L^p(\bRn)}  \leq \|\bnu\|_{[L^q(\bRn\times\bRn)]^k} \|\bkappa\|_{[L^p(\bRn\times\bRn\times\bRn)]^k}
\end{equation}
which completes the proof for \eqref{drega}.

For \eqref{dregb}, let $u\in L^q(\bRn)$. Then, using Minkowski's integral inequality again, we have
$$
\begin{aligned}
 \|\CD^* u \|_{[L^p(\bRn\times\bRn)]^k} &= \Big(\int_{\bRn}\int_{\bRn} \Big|\int_{\bRn} u(\bz)\bkappa(\bz,\bx,\by) d\bz\Big|^p d\by d\bx\Big)^{1/p}\\
                                    &\leq \int_{\bRn} \Big(\int_{\bRn}\int_{\bRn}|u(\bz)\bkappa(\bz,\bx,\by)|^p d\by d\bx\Big)^{1/p} d\bz \\
                                    &\leq \int_{\bRn} |u(\bz)| \Big(\int_{\bRn}\int_{\bRn}|\bkappa(\bz,\bx,\by)|^p\Big)^{1/p} d\bz.
\end{aligned}
$$
Again, applying the H\"{o}lder inequality to the last inequality, we have 
\begin{equation}\label{dregd}
\|\CD^* u\|_{[L^p(\bRn\times\bRn)]^k} \leq \|u\|_{L^q(\bRn)} \|\bkappa\|_{[L^p(\bRn\times\bRn\times\bRn)]^k}
\end{equation}
which completes the proof for \eqref{dregb}. 

The facts that $\CD$ and $\CD^*$ are bounded operators on $L^2(\bRn\times\bRn)^k$ and $L^2(\bRn)$, respectively, follow easily from \eqref{dregc} and \eqref{dregd}, respectively.
\hfill$\Box$

\subsection{Other nonlocal operators}

Other nonlocal operators that mimic the operators of the classical differential vector calculus can be defined. 

\subsubsection{Nonlocal gradient and curl operators}

A nonlocal gradient operator can be defined in a manner similar to Definition \ref{ddiv} for the nonlocal divergence operator.

\begin{defn}[\textbf{Nonlocal gradient operator}]
The action of the {\em nonlocal gradient operator} $\CG: C_c^{\infty}(\bRn\times\bRn)\to [D'(\bRn)]^k$ on any scalar-valued function $\eta\in C_c^{\infty}(\bRn\times\bRn)$ is given by 
\begin{equation}\label{adjgras}
 (\CG\eta)(\bx) = \int_{\bRn}\int_{\bRn}\eta(\by,\bz)\bkappa(\bx,\by,\bz)\, d\bz d\by \qquad\mbox{for a.e. $\bx\in \bRn$}.
\end{equation}
\end{defn}

\begin{prop}
Corresponding to the nonlocal gradient operator \\$\CG: C_c^{\infty}(\bRn\times\bRn)\to [D'(\bRn)]^k$, we have the {\em adjoint operator}\\
 $\CG^\ast: [C_c^{\infty}(\bRn)]^k \to  D'(\bRn\times\bRn)$ whose action on any vector-valued function $\bv\in [C_c^{\infty}(\bRn)]^k$ is given by 
\begin{equation}\label{adjgra}
 (\CG^\ast\bv)(\bx,\by) = \int_{\bRn}\bv(\bz)\cdot\bkappa(\bz,\bx,\by) \,d\bz \qquad\mbox{for a.e. $\bx,\by\in \bRn$}.
\end{equation}
\end{prop}

{\em Proof.}
By definition, the adjoint operator $\CG^\ast$ satisfies
$$
<\bv,\CG\eta>_{[C_c^{\infty}(\bRn)]^k,D'(\bRn\times\bRn)} = <\eta,\CG^* \bv>_{C_c^{\infty}(\bRn\times\bRn),[D'(\bRn)]^k}
$$
for all $\eta\in C_c^{\infty}(\bRn\times\bRn)$ and $\bv\in[ C_c^{\infty}(\bRn)]^k$. Then, the proof of \eqref{adjgra} follows along the same lines of the proof of Proposition \eqref{propadj}.\hfill$\Box$

{\em Remark.}
With $\CG^\ast$ being the adjoint of the nonlocal gradient operator $\CG$, one can identify $-\CG^\ast$ as a nonlocal divergence operator.\hfill$\Box$

{\em Remark.}
We now have the two nonlocal divergence operators $\CD$ and $-\CG^\ast$ and the two nonlocal gradient operators $\CG$ and $-\CD^\ast$. It is natural to have such pairs because of the two types of functions that are needed to describe nonlocality, i.e., functions of two points such as $\bnu(\bx,\by)$ and $\eta(\bx,\by)$ and functions of one point such as $\bv(\bx)$ and $u(\bx)$. Thus, we have the nonlocal divergence and gradient operators $\CD$ and $\CG$ acting on functions of two points and the nonlocal divergence and gradient operators $-\CG^\ast$ and $-\CD^\ast$ acting on functions of one point.\hfill$\Box$

Similar to Proposition \ref{propreg}, one can prove the following proposition.

\begin{prop}[\bf Regularity of $\CG$ and $\CG^*$]\label{propreg2}
Let $1\leq p \leq \infty$ and $1\leq q \leq \infty$ with $\frac{1}{p} + \frac{1}{q} = 1$ and assume that $\bkappa\in [L^p(\bRn\times\bRn\times\bRn)]^k$. Then,
$$
\begin{aligned}
\CG&: \eta\in L^q(\bRn\times\bRn) \longmapsto \CG\eta\in [L^p(\bRn)]^k \\
\CG^*&: \bu\in [L^q(\bRn)]^k \longmapsto \CG^*\bu \in L^p(\bRn\times\bRn) .
\end{aligned}
$$
In particular, 
$$
L^2(\bRn\times\bRn)\xmapsto{\CG} [L^2(\bRn)]^k \xmapsto{\CG^*} L^2(\bRn\times\bRn).
$$
Moreover, $\CG$ and $\CG^*$ are bounded operators on $L^2(\bRn\times\bRn)$ and $[L^2(\bRn)]^k$, respectively.\hfill$\Box$
\end{prop}

In the sequel, we will need the following averaging operator.

\begin{defn}\label{dgbar}
The action of the {\em nonlocal averaging operator}\\
 $\overline{\CG^*} : [C_c^{\infty}(\bRn)]^k \longrightarrow C_c^{\infty}(\bRn)$ on a vector-valued function $\bu\in  [C_c^{\infty}(\bRn)]^k$ is given by
\begin{equation}\label{gbar}
 (\overline{\CG^*}\bu)(\bx) = \int_{\bRn}(\CG^* \bu)(\bx,\bz)\, d\bz.
\end{equation}
\end{defn}

{\em Remark.} A nonlocal curl operator $\CC: [C_c(\bbR^3\times\bbR^3)]^3\to [D'(\bbR^3)]^3$ is given by its action action on any vector-valued function $\bnu\in [C_c^{\infty}(\bbR^3\times\bbR^3)]^3$ as
$$
 (\CC\betab)(\bx) = \int_{\bbR^3}\int_{\bbR^3}\betab(\by,\bz)\times\bkappa(\bx,\by,\bz)\, d\bz d\by \qquad\mbox{for a.e. $\bx\in \bbR^3$}.
$$
The corresponding nonlocal adjoint operator $\CC^\ast: [C_c^{\infty}(\bbR^3)]^3\to [D'(\bbR^3\times\bbR^3)]^3$, which is also a nonlocal curl operator, is given by its action on any vector-valued function $\bu\in [C_c^{\infty}(\bbR^3)]^3$ as
$$
 (\CC^\ast\bv)(\bx,\by) = \int_{\bbR^3}\bkappa(\bz,\bx,\by)\times\bu(\bz) \,d\bz \qquad\mbox{for a.e. $\bx,\by\in \bbR^3$}.
$$
Regularity results similar to those proved in Propositions\eqref{propreg} and \eqref{propreg2} for the nonlocal divergence and gradient operators hold for the nonlocal curl operator $\CC$.\hfill$\Box$

\subsubsection{Nonlocal divergence of a tensor and gradient of a vector}

The nonlocal divergence operator $\CD$ can also be applied to a tensor-valued function yielding a vector-valued function.
 
\begin{defn} [\textbf{Nonlocal divergence of a tensor}]
The action of the nonlocal divergence operator 
$\CD:[C_c^{\infty}(\bRn\times\bRn)]^{k\times k}\longrightarrow[D'(\bRn)]^k$ on the tensor-valued function $\bPsi \in [C_c^{\infty}(\bRn\times\bRn)]^{k\times k}$ is defined by 
\begin{equation}\label{divten}
 (\CD\bPsi)(\bx) = \int_{\bRn}\int_{\bRn}\bPsi(\by,\bz)\bkappa(\bx,\by,\bz)\,d\bz d\by
 \qquad\mbox{for a.e. $\bx\in \bRn$}.
\end{equation}
\end{defn}
Here, $\bPsi\bkappa$ represents a matrix-vector product. The components of the vector $\CD\bPsi$ are the nonlocal divergences of the corresponding rows of $\bPsi$.

\begin{prop}
The action of the
nonlocal adjoint operator\\ $\CD^*:[C_c^{\infty}(\bRn)]^k\longrightarrow[D'(\bRn\times\bRn)]^{k\times k}$ on the vector-valued function $\bu\in [C_c^{\infty}(\bRn)]^k$ is given by 
\begin{equation}\label{adjten}
 (\CD^*\bu) (\bx,\by) = \int_{\bRn} \bu(\bz)\otimes\bkappa(\bz,\bx,\by) d\bz \qquad\mbox{for a.e. $\bx,\by\in \bRn$}.
\end{equation}
\end{prop}

{\em Proof.}
By the definition of adjoint operator, we have 
$$
< \bu,\CD\bPsi>_{[C_c^{\infty}(\bRn)]^k,[D'(\bRn)]^k} = <\bPsi,\, \CD^*\bu>_{[C_c^{\infty}(\bRn\times\bRn)]^{k\times k},[D'(\bRn\times\bRn)]^{k\times k}} 
$$
for all $\bPsi\in C_c^{\infty}[\bRn\times\bRn)]^{k\times k}$ and $\bu\in[C_c^{\infty}(\bRn)]^k$.
This is equivalent to 
\begin{equation}\label{adjten2}
\begin{aligned}
\int_{\bRn} \Big(\int_{\bRn}\int_{\bRn}\bPsi(\by,\bz)&\bkappa(\bx,\by,\bz)\,d\bz d\by \Big)\cdot \bu(\bx)\, d\bx 
\\&= 
\int_{\bRn}\int_{\bRn} \bPsi(\bx,\by) : (\CD^*\bu)(\bx,\by) \,d\by d\bx.
\end{aligned}
\end{equation}
After switching $\bx$ and $\bz$ and then $\bx$ and $\by$ in the left-hand side, we have
$$
\begin{aligned}
\int_{\bRn} \int_{\bRn}\int_{\bRn}&\big(\bPsi(\by,\bz)\bkappa(\bx,\by,\bz)\big) \cdot \bu(\bx)\,d\bz d\by d\bx 
\\&=
\int_{\bRn}\int_{\bRn}\int_{\bRn}\big(\bPsi(\by,\bx)\bkappa(\bz,\by,\bx)\big) \cdot \bu(\bz)\,d\bx d\by d\bz
\\&=
\int_{\bRn}\int_{\bRn}\int_{\bRn}\big(\bPsi(\bx,\by)\bkappa(\bz,\bx,\by)\big) \cdot\bu(\bz)\,d\by d\bx d\bz
\\&=
\int_{\bRn}\int_{\bRn}\bPsi(\bx,\by):\Big(\int_{\bRn}
\bu(\bz) \otimes\bkappa(\bz,\bx,\by)\,d\bz\Big)\,d\by d\bx, 
\end{aligned}
$$
where, for the last equality, we rearranged the tensor-vector products. Then, because $\bPsi\in [C_c^{\infty}(\bRn\times\bRn)]^{k\times k}$ is arbitrary, we obtain \eqref{adjten} from \eqref{adjten2}.\hfill$\Box$

The nonlocal gradient operator $\CG$ can also be applied to a vector-valued function yielding a tensor-valued function.

\begin{defn} [\textbf{Nonlocal gradient of a vector}]
The action of the nonlocal gradient operator 
$\CG:[C_c^{\infty}(\bRn\times\bRn)]^{k}\longrightarrow[D'(\bRn)]^{k\times k}$ on the vector-valued function $\bnu \in [C_c^{\infty}(\bRn\times\bRn)]^{k}$ is defined by 
\begin{equation}\label{gradvec}
 (\CG\bnu)(\bx) = \int_{\bRn}\int_{\bRn}\bnu(\by,\bz)\otimes\bkappa(\bx,\by,\bz)\,d\bz d\by
 \qquad\mbox{for a.e. $\bx\in \bRn$}.
\end{equation}
\end{defn}
We then have that the action of the nonlocal adjoint operator $\CG^*:[C_c^{\infty}(\bRn)]^{k\times k}\rightarrow[D'(\bRn\times\bRn)]^{k}$ on the tensor-valued function\\
 $\BU\in [C_c^{\infty}(\bRn)]^{k\times k}$ is given by 
\begin{equation}\label{adjten2}
 (\CG^*\BU) (\bx,\by) = \int_{\bRn} \BU(\bz)\bkappa(\bz,\bx,\by) d\bz \qquad\mbox{for a.e. $\bx,\by\in \bRn$}.
\end{equation}

\subsubsection{Nonlocal Laplacian operators}

With $\CD$ and $-\CD^\ast$ denoting nonlocal divergence and gradient operators, respectively, their composition $-\CD\CD^\ast$ can be viewed as a nonlocal Laplacian operator. 
The following proposition provides the explicit form of this operator.

\begin{prop}
The nonlocal Laplacian operator of a scalar-valued function $u(\bx)$ is given by
\begin{equation}\label{nonlapa}
   -(\CD\CD^\ast u)(\bx)=
  -\int_{\bRn}\int_{\bRn}\int_{\bRn}
  u(\bw)\bkappa(\bw,\by,\bz)\cdot\bkappa(\bx,\by,\bz)\, d\bw d\bz d\by
\end{equation}
whereas the nonlocal Laplacian operator of a vector-valued function $\bu(\bx)$ is given by
\begin{equation}\label{nonlapb}
-\CD(\CD^\ast {\bf u}) =
-\int_{\bRn}\int_{\bRn}\int_{\bRn}\bu(\bw)\bkappa(\bw,\by,\bz)
\cdot\bkappa(\bx,\by,\bz)\, d\bw d\bz d\by
\end{equation}

\end{prop}

{\em Proof.}
From \eqref{divop} and \eqref{adjop}, we have that
$$
\begin{aligned}
\CD\CD^\ast u &= \int_{\bRn}\int_{\bRn} (\CD^\ast u)(\by,\bz)
\cdot\bkappa(\bx,\by,\bz)\,d\bz d\by
\\&=
\int_{\bRn}\int_{\bRn} \Big[
\int_{\bRn} u(\bw)\bkappa(\bw,\by,\bz)\, d\bw
\Big]\cdot\bkappa(\bx,\by,\bz)\,d\bz d\by.
\end{aligned}
$$
In the same manner, \eqref{nonlapb} follows from \eqref{divten} and  \eqref{adjten}.
\hfill$\Box$

\subsubsection{Identities of the nonlocal calculus}\label{sec:ident}

We begin with some identities that mimic those of the classical vector calculus. The first set of identities do not require any further conditions on the divergence kernel $\bkappa(\bx,\by,\bz)$.
\begin{subequations}\label{ident1}
\begin{prop}\label{propident1}

{\rm(i)} For $u(\bx)=a$ and $\bu(\bx)={\bf a}$, where $a$ and ${\bf a}$ are scalar and vector constants, respectively, we have
\begin{equation}\label{ident1a}
(\CD^\ast a)(\bx,\by) ={\bf 0},\qquad
(\CG^\ast{\bf a})(\bx,\by) =0,\qquad\mbox{and}\qquad
(\CC^\ast{\bf a})(\bx,\by) ={\bf 0}.
\end{equation}

{\rm(ii)} For the vector-valued functions $\bu(\bx)$ and $\bnu(\bx,\by)$, we have 
\begin{equation}\label{ident1b}
\begin{aligned}
   (\CD\CD^\ast {\bf u})(\bx)&= (\CC \CC^\ast {\bf u})(\bx)+
   (\CG \CG^\ast {\bf u})(\bx)
\\
   (\CG^\ast\CG\bnu)(\bx,\by) 
&= (\CC^\ast\CC\bnu)(\bx,\by)+(\CD^\ast \CD\bnu)(\bx,\by).
\end{aligned}
\end{equation}

{\rm(iii)} For the vector-valued functions $\bu(\bx)$ and $\bnu(\bx,\by)$, we have 
\begin{equation}\label{ident1c}
 \CD\bnu = \mbox{\em trace} (\CG\bnu) 
 \qquad\mbox{and}\qquad
  \CG^\ast \bu = \mbox{\em trace} (\CD^\ast \bu).
\end{equation}
\end{prop}
\end{subequations}
{\em Proof.} (i) Using \eqref{divthm4}, we have that
$$
   (\CD^\ast a)(\bx,\by) = \int_\bRn a\bkappa(\bz,\bx,\by)\,d\bz=
   a\int_\bRn \bkappa(\bz,\bx,\by)\,d\bz={\bf0}
   \qquad\forall\, a.
$$
so that $(\CD^\ast a)(\bx,\by) ={\bf0}$. The other two results in \eqref{ident1a} are proved in a similar manner.

(ii)
We have that
$$
\begin{aligned}
\CD(\CD^\ast {\bf u})- \CG(\CG^\ast {\bf u})&
=\int_{\bbR^3}\int_{\bbR^3}(\CD^\ast {\bf u})(\by,\bz)\bkappa(\bx,\by,\bz)\,d\bz d\by 
\\& \qquad-  \int_{\bbR^3}\int_{\bbR^3}(\CG^\ast {\bf u})(\by,\bz) \bkappa(\bx,\by,\bz)\,d\bz d\by 
\\& =
\int_{\bbR^3}\int_{\bbR^3}\int_{\bbR^3}\Big( \bu(\bw)\otimes\bkappa(\bw,\by,\bz)\Big)
\bkappa(\bx,\by,\bz)\, d\bw d\bz d\by 
\\& \qquad- 
\int_{\bbR^3}\int_{\bbR^3}\int_{\bbR^3} \Big(\bu(\bw)\cdot\bkappa(\bw,\by,\bz)\Big)
\bkappa(\bx,\by,\bz)\,d\bw d\bz d\by
\\&=
\int_{\bbR^3}\int_{\bbR^3}\int_{\bbR^3}\Big[
\bu(\bw)\Big(\bkappa(\bw,\by,\bz)
\cdot\bkappa(\bx,\by,\bz)\Big)
\\& \qquad
-\Big(\bu(\bw)\cdot\bkappa(\bw,\by,\bz)\Big)
\bkappa(\bx,\by,\bz)\Big]\,d\bw d\bz d\by
\\&
=\int_{\bbR^3}\int_{\bbR^3} \int_{\bbR^3}
\bkappa(\bw,\by,\bz)\times\bu(\bw)
\times \bkappa(\bx,\by,\bz)\,d\bw d\bz d\by  
\\&
=\int_{\bbR^3}\int_{\bbR^3} \Big[\int_{\bbR^3}
\bkappa(\bw,\by,\bz)\times\bu(\bw)
\,d\bw\Big] \times \bkappa(\bx,\by,\bz)\,d\bz d\by
\\&
=\int_{\bbR^3}\int_{\bbR^3} (\CC^\ast \bu)(\by,\bz) \times \bkappa(\bx,\by,\bz)\,d\bz d\by
=\CC(\CC^\ast\bu) , 
\end{aligned}
$$
where the first two equalities follow from the definitions of the operators $\CD$, $\CD^\ast$, $\CG$, and $\CG^\ast$, the third and fourth equalities follow from the standard vector identities $({\bf a}\otimes{\bf b})\cdot{\bf c} = {\bf a}({\bf b}\cdot{\bf c})$ and ${\bf a}\times({\bf b}\times{\bf c}) = {\bf b}({\bf a}\cdot{\bf c}) - {\bf c}({\bf a}\cdot{\bf b})$, respectively, the fifth equality is a tautology, and the last two inequalities follow from the definitions of the operators $\CC^\ast$ and $\CC$. 
The second identity in \eqref{ident1b} is proved in a similar manner.


(iii)
The proofs of the identities in \eqref{ident1c} follow easily from the definitions of the operators and of the matrix trace, e.g.,
$$
\begin{aligned}
\mbox{trace} (\CG\bnu) &= \mbox{trace}\Big(\int_{\bRn}\int_{\bRn} \bnu(\by,\bz)\otimes \bkappa(\bx,\by,\bz)\,d\bz d\by\Big)
\\&= \int_{\bRn}\int_{\bRn} \mbox{trace}\Big(\bnu(\by,\bz)\otimes \bkappa(\bx,\by,\bz)\Big)\,d\bz d\by
\\&
 =  \int_{\bRn}\int_{\bRn}  \bnu(\by,\bz)\cdot \bkappa(\bx,\by,\bz)\,d\bz d\by =\CD\bnu
\end{aligned}
$$
with the proof of the second identity in \eqref{ident1c} following in a similar manner.~\hfill$\Box$

Unlike the identities \eqref{ident1}, the second set of identities do require additional conditions on the divergence kernel.

\begin{subequations}\label{ident2}
\begin{prop}\label{propident2}

{\rm(i)} For $u(\bx)=a$ and $\bu(\bx)={\bf a}$, where $a$ and ${\bf a}$ are scalar and vector constants, respectively, we have
\begin{equation}\label{ident2a}
(\CD{\bf a})(\bx) =0,\qquad
(\CG a)(\bx) ={\bf 0},\qquad\mbox{and}\qquad
(\CC{\bf a})(\bx) ={\bf 0}
\end{equation}
if and only if
\begin{equation}\label{ident2b}
\int_\bRn\int_\bRn \bkappa(\bx,\by,\bz)\,d\bz d\by={\bf 0}
\qquad
\forall\,\bx.
\end{equation}

{\rm(ii)} For any $v(\bx)$ and $\bu(\bx)$, we have that
\begin{equation}\label{ident2c}
\CD\big(\CC^\ast\bu\big)(\bx)=0\qquad\mbox{and}\qquad
\CC\big(\CD^\ast v\big)(\bx)={\bf0}
\end{equation}
if and only if
\begin{equation}\label{ident2d}
\int_\bRn\int_\bRn
 \bkappa(\bx,\by,\bz)\times\bkappa(\bw,\by,\bz)
\,d\bz d\by={\bf 0}
\qquad\forall\,\bx,\,\bw.
\end{equation}

{\rm(iii)} For any $\bnu(\bx,\by)$ and $\eta(\bx,\by)$ we have that
\begin{equation}\label{ident2e}
\CG^\ast \big(\CC\bnu\big)(\bx,\by)=0\qquad\mbox{and}\qquad
\CC^\ast\big(\CG\eta\big)(\bx,\by)=0
\end{equation}
if and only if
\begin{equation}\label{ident2f}
\int_\bRn
\bkappa(\bz,\bw,\br)\times\bkappa(\bz,\bx,\by)
\,d\bz = {\bf0}
\qquad\forall\,\bx,\,\by,\,\bw,\,\br.
\end{equation}
\end{prop}
\end{subequations}

{\em Proof.}
(i) For $\CD{\bf a}$, we have
$$
(\CD{\bf a})(\bx) =
\int_\bRn\int_\bRn {\bf a}\cdot\bkappa(\bx,\by,\bz)\,d\bz d\bz
={\bf a}\cdot\int_\bRn\int_\bRn \bkappa(\bx,\by,\bz)\,d\bz d\by
 \qquad\forall\, {\bf a}
$$
so that $(\CD{\bf a})(\bx) = 0$ if and only if \eqref{ident2b} holds. The other two results in \eqref{ident2a} are proved in a similar manner.

{\rm(ii)} From the definitions of the operators $\CC$ and $\CD^\ast$, we have that
$$
\begin{aligned}
\CD\big(\CC^\ast \bu\big)(\bx)
&=\int_{\bbR^3}\int_{\bbR^3} (\CC^\ast \bu)(\by,\bz)\cdot
\bkappa(\bx,\by,\bz)\,d\bz d\by
\\&
=\int_{\bbR^3}\int_{\bbR^3} \Big[
\int_{\bbR^3} \bkappa(\bw,\by,\bz) \times{\bf u}(\bw) \,d\bw
\Big]\cdot
\bkappa(\bx,\by,\bz)\,d\bz d\by
\\&
=\int_{\bbR^3} \bu(\bw)\cdot\Big[ \int_{\bbR^3}\int_{\bbR^3}
 \bkappa(\bw,\by,\bz)\times
\bkappa(\bx,\by,\bz)\,d\bz d\by\Big]\,d\bw.
\end{aligned}
$$
Because $\bu(\bx)$ is arbitrary, the first result in \eqref{ident2c} follows; the second results follows in a similar manner.

{\rm(iii)} From the definitions of the operators $\CG^\ast$ and $\CC$, we have that 
$$
\begin{aligned}
\CG^\ast\big(\CC\bnu\big)(\bx,\by)
&=\int_{\bbR^3}(\CC\bnu)(\bz)\cdot\bkappa(\bz,\bx,\by)\,d\bz
\\&=
\int_{\bbR^3}\Big[
\int_{\bbR^3}\int_{\bbR^3} 
\bnu(\bw,\br)\times\bkappa(\bz,\bw,\br) \,d\bw d\br
\Big]\cdot\bkappa(\bz,\bx,\by)\,d\bz
\\&=
\int_{\bbR^3}\int_{\bbR^3} \bnu(\bw,\br)\cdot\Big[\int_{\bbR^3}
\bkappa(\bz,\bw,\br)\times\bkappa(\bz,\bx,\by)
\,d\bz\Big]\,d\bw d\br.
\end{aligned}
$$
Because $\bnu(\bx,\by)$ is arbitrary, the first result in \eqref{ident2e} follows; the second results follows in a similar manner. \hfill$\Box$

\subsubsection{Theorems of the nonlocal calculus}

We next consider the nonlocal analog of the divergence theorem of the classical vector calculus.

\begin{thm}\label{divthmthm}
Let $\Omega\subseteq\bRn$. Then,
\begin{equation}\label{divthmnld}
\int_\Omega (\CD\bnu)(\bx)\,d\bx
= - \int_{\bRn\setminus\Omega} (\CD\bnu)(\bx)\,d\bx.  
\end{equation} 
\end{thm}

{\em Proof.} The proof of \eqref{divthmnld} is basically a tautology because we defined the nonlocal operator $\CD$ so that it satisfies a nonlocal divergence theorem. In fact, \eqref{divthmnld} follows easily from \eqref{divthm5}.\hfill$\Box$

{\em Remark.} The integral of the classical local divergence of a vector over and arbitrary domain $\Omega$ is equal to the flux of that vector out of $\Omega$ which is given by an integral over the boundary of $\Omega$ of the normal component of the vector. Nonlocality results in the flux out of $\Omega$ to be given by a volume integral over the complement of $\Omega$ as is indicated in \eqref{divthmnld}.\hfill$\Box$

{\em Remark.} Analogous theorems hold for the operators $\CG$ and $\CC$, i.e., $\int_\Omega (\CG\eta)(\bx)\,d\bx = - \int_{\bRn\setminus\Omega} (\CG\eta)(\bx)\,d\bx$ and $\int_\Omega (\CC\betab)(\bx)\,d\bx = - \int_{\bbR^3 \setminus\Omega} (\CD\betab)(\bx)\,d\bx$.\hfill$\Box$

Finally, we derive the nonlocal Green's identities which again mimic the classical Green's identities of the classical vector calculus. We begin with an integration by parts formula.

\begin{lem}\label{greenide}
Given any functions $u(\bx)$ and $\bnu(\bx,\by)$, we have that
\begin{equation}\label{greenibp}
\int_\bRn  u(\bx) \CD(\bnu)(\bx)\,d\bx
 -
 \int_\bRn\int_\bRn \bnu(\bx,\by)\cdot (\CD^\ast u)(\bx,\by)
 \, d\by d\bx
 =0.
\end{equation}
\end{lem}

{\em Proof.}
We have that
$$
\begin{aligned}
  0 &= \int_\bRn\int_\bRn\int_\bRn u(\bx)
 \bnu(\by,\bz)\cdot\bkappa(\bx,\by,\bz)\,d\bz d\by d\bx
\\&\qquad - \int_\bRn\int_\bRn\int_\bRn
 u(\bx)\bnu(\by,\bz)\cdot\bkappa(\bx,\by,\bz)\,d\bz d\by d\bx
\\&=
\int_\bRn\int_\bRn\int_\bRn u(\bx)
 \bnu(\by,\bz)\cdot\bkappa(\bx,\by,\bz)\,d\bz d\by d\bx
\\&\qquad -
\int_\bRn\int_\bRn\int_\bRn
 u(\bz)\bnu(\bx,\by)\cdot\bkappa(\bz,\bx,\by)\,d\bz d\by d\bx
\\&=
\int_\bRn u(\bx) \Big[\int_\bRn\int_\bRn
\bnu(\by,\bz)\cdot\bkappa(\bx,\by,\bz)\,d\bz d\by\Big]\,d\bx
\\&\qquad -
\int_\bRn\int_\bRn \bnu(\bx,\by)\cdot \Big[\int_\bRn
 u(\bz)\bkappa(\bz,\bx,\by)\,d\bz\Big]\, d\by d\bx
 \\& = 
 \int_\bRn  u(\bx) \CD(\bnu)(\bx)\,d\bx
 -
 \int_\bRn\int_\bRn \bnu(\bx,\by)\cdot (\CD^\ast u)(\bx,\by)
 \, d\by d\bx
\end{aligned}
$$
where the first equality is a tautology, the second equality follows from a cyclic replacement of the integration variables ($\bx$, $\by$, $\bz$ $\to$ $\bz$, $\bx$, $\by$) in the second integral, the third equality is again a tautology, and the last follows from the definition of the operators $\CD$ and $\CD^*$. Thus, \eqref{greenibp} is proven.
\hfill$\Box$

\begin{thm}[\bf Green's identities]\label{greenidenthm}
Given functions $u(\bx)$ and $v(\bx)$, we have the {\em nonlocal Green's first identity}
\begin{subequations}\label{greeniden}
\begin{equation}\label{greeniden1}
\int_\bRn  u \CD(\CD^\ast v)\,d\bx
 -
 \int_\bRn\int_\bRn (\CD^\ast v)\cdot (\CD^\ast u)
 \, d\by d\bx =0
\end{equation}
the {\em nonlocal Green's second identity}
\begin{equation}\label{greeniden2}
\int_\bRn  u \CD(\CD^\ast v)\,d\bx -
 \int_\bRn  v \CD(\CD^\ast u)\,d\bx =0.
\end{equation}
\end{subequations}
\end{thm}

{\em Proof.}
Setting $\bnu(\bx,\by) = (\CD^\ast v)(\bx,\by)$ in \eqref{greenibp} easily results in \eqref{greeniden1}. Then, \eqref{greeniden2} follows by reversing the roles of $u$ and $v$ in \eqref{greeniden1} and then subtracting the result from \eqref{greeniden1}. \hfill$\Box$

{\em Remark.} Analogous theorems hold for the pairs of operators $\CG$ and $\CG^\ast$ and $\CC$ and $\CC^\ast$.

The following results are obvious consequences of \eqref{greeniden1} and \eqref{greeniden2}.

\begin{cor}
Given a subdomain $\Omega\subseteq\bRn$ and functions $u(\bx)$ and $v(\bx)$, we have that
\begin{subequations}\label{greeniden}
\begin{equation}\label{greeniden1aa}
\int_\Omega  u \CD(\CD^\ast v)\,d\bx
 -
 \int_\bRn\int_\bRn (\CD^\ast v)\cdot (\CD^\ast u)
 \, d\by d\bx=
 - \int_{\bRn\setminus\Omega}  u \CD(\CD^\ast v)\,d\bx
\end{equation}
and
\begin{equation}\label{greeniden2bb}
\int_\Omega  u \CD(\CD^\ast v)\,d\bx -
 \int_\Omega  v \CD(\CD^\ast u)\,d\bx =
 - \int_{\bRn\setminus\Omega}  u \CD(\CD^\ast v)\,d\bx
 + \int_{\bRn\setminus\Omega}  v \CD(\CD^\ast u)\,d\bx
 \hfill\Box.
\end{equation}
\end{subequations}
\end{cor}

\section{Special case of the nonlocal operators}
\label{sec_special}

The general forms of the nonlocal divergence operator and its adjoint operator are given in Definition \ref{ddiv} and Proposition \ref{propadj}. Here, we consider a simplified version of these operators which leads to the nonlocal vector calculus of \cite{dglz-nlc} and which has proven to be useful \cite{dglz-nld,dglz-ps}.

The simplification is effected by a special case of the Schwartz kernel given by
\begin{equation}\label{specker}
\bkappa_{\balpha}(\bx,\by,\bz) = \delta(\bx-\bz)\balpha(\bx,\by) + \delta(\bx-\by)\balpha(\bx,\bz)
\end{equation}
for an vector-valued function $\balpha(\bx,\by)\in[L^1(\bRn\times\bRn)]^k$. Here, $\delta(\cdot)$ denotes the Dirac delta function. First, we verify that $\bkappa_\alpha$ satisfies \eqref{divthm4}.

\begin{prop}\label{skdiv}
The specialized Schwartz kernel $\bkappa_{\balpha}$ satisfies \eqref{divthm4}, i.e., 
\begin{equation}\label{divthm4a}
\int_{\bRn}\bkappa_{\balpha}(\bx,\by,\bz)\,d\bx = 0, 
\end{equation}
if and only if $\balpha(\bx,\by)$ is antisymmetric, i.e., if and only if $\balpha(\bx,\by) = -\balpha(\by,\bx)$ for all $\bx$ and $\by$.
\end{prop}

{\em Proof.}
We have 
$$
\begin{aligned}
\int_{\bRn}\bkappa_{\balpha}(\bx,\by,\bz)\,d\bx &= \int_{\bRn}\big(\delta(\bx-\bz)\balpha(\bx,\by) + \delta(\bx-\by)\balpha(\bx,\bz)\big)\,d\bx
\\&= \balpha(\bz,\by) + \balpha(\by,\bz) 
\end{aligned}
$$
so that the result follows.
\hfill$\Box$

{\em Remark.}
In \cite{dglz-nld,dglz-nlc,dglz-ps}, the antisymmetry of $\balpha(\bx,\by)$ is {\em assumed.} Here, we have shown that this condition is necessary and sufficient for the operator $\CD$ to be a nonlocal divergence operator in the sense that \eqref{divthm} (and therefore \eqref{divthm4}) is satisfied.
\hfill$\Box$

\begin{thm}[\textbf{Specialized nonlocal divergence operator and its adjoint}]\label{thmsdiv}
For the specialized kernel $\bkappa_{\balpha}$ given by \eqref{specker}, the action of the nonlocal divergence operator $(\CD_{\balpha}\bnu)(\bx):\Us\rightarrow\Vsp$ on a function $\bnu(\bx,\by)\in \Us$ is given by
\begin{equation}\label{sgrad}
 (\CD_{\balpha}\bnu)(\bx) = \int_{\bRn}\big(\bnu(\bx,\by) + \bnu(\by,\bx)\big)\cdot\balpha(\bx,\by)\,d\by.
\end{equation}
Moreover, the action of the adjoint operator $(\CD^\ast_{\balpha} u)(\bx,\by): \Vs \rightarrow \Usp$ on a function $u(\bx)\in \Vs$ is given by 
\begin{equation}\label{sgradadj}
 (\CD_{\balpha}^\ast u)(\bx,\by) = -\big(u(\by) - u(\bx)\big)\balpha(\bx,\by) .
\end{equation}
\end{thm}

{\em Proof.}
Setting $\bkappa=\bkappa_{\balpha}$ in \eqref{divop}, we have
$$
\begin{aligned}
 \int_{\bRn}\int_{\bRn}\bkappa_{\balpha}&(\bx,\by,\bz)\cdot\nu(\by,\bz)\,d\bz d\by \\
&= \int_{\bRn}\int_{\bRn}\Big(\delta(\bx-\bz)\balpha(\bx,\by) + \delta(\bx-\by)\balpha(\bx,\bz)\Big)\cdot\bnu(\by,\bz)\, d\bz d\by \\
&= \int_{\bRn}\balpha(\bx,\by)\cdot\bnu(\by,\bx)\,d\by + \int_{\bRn}\balpha(\bx,\bz)\cdot\bnu(\bx,\bz)\,d\bz\\
&=\int_{\bRn}\big(\bnu(\bx,\by)+  \bnu(\by,\bx)\big) \cdot\balpha(\bx,\by)    \,d\by,
\end{aligned}
$$
where the last equality follows by replacing the integration variable from $\bz$ to $\by$ in the second integral preceding the equality. Thus, we have \eqref{sgrad}.

Setting $\bkappa=\bkappa_{\balpha}$ in \eqref{adjop}, we have 
$$
\begin{aligned}
\int_{\bRn}u(\bz)\bkappa_{\balpha}(\bz,\bx,\by)\,d\bz
                   &= \int_{\bRn}u(\bz)\Big(\delta(\bz-\by)\balpha(\bz,\bx) + \delta(\bz-\bx)\balpha(\bx,\by)\Big)d\bz
              \\     &= u(\by)\balpha(\by,\bx) + u(\bx)\balpha(\bx,\by) =
              -\big(u(\by) - u(\bx)\big)\balpha(\bx,\by) ,
\end{aligned}
$$
completing the proof of \eqref{sgradadj}.\hfill$\Box$

{\em Remark.}
For the kernel $\bkappa_{\balpha}$, we have  
\begin{equation}\label{mmmmmm}
\begin{aligned}
(\CD_{\balpha}\bPsi)(\bx) &= \int_{\bRn}\big(\bPsi(\by,\bx) + \bPsi(\bx,\by)\big) \balpha(\bx,\by)\,d\by
\\[1ex]
(\CD^*_{\balpha}\bu)(\bx,\by) &= -\big(\bu(\by) - \bu(\bx)\big)\otimes\balpha(\bx,\by)
\\[1ex]
 (\CG_{\balpha}\eta)(\bx) &= \int_{\bRn}\big(\eta(\by,\bx) + \eta(\bx,\by)\big)\balpha(\bx,\by)\, d\by
\\[1ex]
 (\CG_{\balpha}^* \bu)(\bx,\by) & = -\big(\bu(\by)-\bu(\bx)\big)\cdot\balpha(\bx,\by) = \mbox{trace}\, (\CD^*_{\balpha} \bu)
\\[1ex]
 (\overline{\CG^*_{\balpha}}\bu)(\bx) &= -\int_{\bRn}\big(\bu(\bz) - \bu(\bx)\big) \cdot \balpha(\bx,\bz) \,d\bz
\\[1ex]
 (\CC_{\balpha}\bnu)(\bx) &= \int_{{\bbR^3}}\balpha(\bx,\by)\times\big(\bnu(\by,\bx) + \bnu(\bx,\by)\big)\, d\by
\\[1ex]
 (\CC_{\balpha}^* \bu)(\bx,\by) & = -\big(\bu(\by)-\bu(\bx)\big)\times\balpha(\bx,\by) 
\\[1ex]
  (\CD_{\balpha}\CD_{\balpha}^\ast u)(\bx) &=-2
  \int_{\bRn} \big(u(\by) - u(\bx)\big)
  \balpha(\bx,\by)\cdot \balpha(\bx,\by)\,d\by
\end{aligned}
\end{equation}
for a tensor-valued $\bPsi(\bx,\by)$, a vector valued function $\bu(\bx)$, and a scalar-valued functions $\eta(\bx,\by)$.
\hfill$\Box$

We now want to examine the identities of Section \ref{sec:ident}  in the context of the specialized kernel $\bkappa_{\balpha}(\bx,\by,\bz)$. Instead of verifying the assumptions of Proposition \ref{propident2}, we directly examine those identities for the kernel \eqref{specker}. Of course, because of Propositions \ref{propident1} and \ref{skdiv}, we have all the identities \eqref{ident1} hold for the operators $\CD_{\balpha}$, $\CD_{\balpha}^\ast$, $\CG_{\balpha}$, etc. Thus, we need only address the identities \eqref{ident2}.

From the definitions of the relevant operators and the antisymmetry of $\balpha(\bx,\by)$, we have that
$$
\begin{aligned}
(&\CD_{\balpha} \CC_{\balpha}^\ast\bu)(\bx)
\\&=\int_{{\bbR^3}}
\Big(\balpha(\bx,\by) \times \big( \bu(\bx) - \bu(\by) \big)
+\balpha(\by,\bx) \times \big( \bu(\by) - \bu(\bx) \big)\Big)
 \cdot \balpha(\bx,\by)\,d\by
\\&=
2\int_{{\bbR^3}} 
\Big(\balpha(\bx,\by) \times \big( \bu(\bx) - \bu(\by) \big)
\Big)
 \cdot \balpha(\bx,\by)
\,d\by = {\bf0}.
\end{aligned}
$$
Similarly, one can show that $(\CC_{\balpha} \CD_{\balpha}^\ast\bu)(\bx)={\bf 0}$. Also, we have that
$$
\begin{aligned}
(\CG_{\balpha}^\ast &\CC_{\balpha}\bnu)(\bx,\by)
\\&=-\int_{{\bbR^3}}\Big(
\balpha(\bx,\by)\times\big(\bnu(\bx,\by)+\bnu(\by,\bx)\big)
\\&
\qquad\qquad-
\balpha(\by,\bx)\times\big(\bnu(\by,\bx)+\bnu(\bx,\by)\big)
\Big)\cdot\balpha(\bx,\by)\,d\by
\\&=
-2\int_{{\bbR^3}}\Big(
\balpha(\bx,\by)\times\big(\bnu(\bx,\by)+\bnu(\by,\bx)\big)
\Big)\cdot\balpha(\bx,\by)\,d\by={\bf0}
\end{aligned}
$$
and similarly that $(\CC^\ast\CG\eta)(\bx,\by) = {\bf 0}$. Finally, we have that, in general,
$$
  (\CD{\bf a})(\bx)=
  2{\bf a}\cdot\int_{\bRn} \balpha(\bx,\by)\,d\by \ne0
$$
and likewise $(\CG a)(\bx)\ne0$ and $(\CC{\bf a})(\bx)\ne{\bf 0}$. The following proposition summarizes these results.

\begin{prop}
In general, for the specialized Schwartz kernel $\bkappa_{\balpha}$, we have that 
\eqref{ident1a}, \eqref{ident1b}, \eqref{ident1c}, \eqref{ident2c}, and \eqref{ident2e} hold, but \eqref{ident2a} does not hold. \hfill$\Box$
\end{prop}

\section{The peridynamics model for solid mechanics}
\label{sec_peridynamics}

We consider the state-based peridynamics model introduced in \cite{Silling2007} for the dynamics of isotropic heterogeneous solids. To simplify the presentation, we provide a direct description of this model without adhering to the notation used in \cite{Silling2007}. Our goal is to show how the peridynamics model can be expressed in terms of the nonlocal operators introduced and discussed in Section \ref{gennlc}.

\subsection{The peridynamics model for solid mechanics}\label{peridyn}

Let $\Omega$ denote a domain in $\bbR^3$, $\bu(\bx,t)$ the displacement vector field, $\rho(\bx)$ the mass density, and $b(\bx,t)$ a prescribed body force density. Let $B_\varepsilon(\bx)$ denote the ball centered at $\bx$ having radius $\varepsilon$; here, $\varepsilon$ denotes the peridynamics horizon. Then, the peridynamic equation of motion is given by
\begin{equation}\label{pd1}
 \rho(\bx) \ddot{\bu}(\bx,t) = \int_{B_\varepsilon(\bx)}\big( \TTT(\bx,\by-\bx) - \TTT(\by, \bx - \by)\big) \,d\by + \bb(\bx,t),
\end{equation}
where 
\begin{equation}\label{pd2}
  \TTT(\bx,\by-\bx) = \ttt(\bx,\by)\bgamma(\bx,\by)
\end{equation}
with
\begin{equation}\label{pd3}
 \bgamma(\bx,\by) = \frac{\bu(\by) + \by - (\bu(\bx) + \bx)}{|\bu(\by) + \by - (\bu(\bx) + \bx)|}
\end{equation}
and 
\begin{equation}\label{pd4}
\begin{aligned}
  \ttt&(\bx,\by) = \frac{3k(\bx)}{m} w(|\by-\bx|) |\by-\bx|  \theta(\bx) \\
             &+ \frac{15\mu(\bx)}{m} w(|\by - \bx|)\Big(|\bu(\by) + \by - (\bu(\bx) + \bx)| - |\by-\bx| -\frac{1}{3}|\by-\bx| \theta(\bx) \Big).
\end{aligned}
\end{equation}
In \eqref{pd1}--\eqref{pd4}, $\bgamma$ represents the direction of the force density that the particle at position $\by$ exerts on the particle at position $\bx$ and $\ttt$ represents the magnitude of that force density. The first term in $\ttt$ is the \emph{hydrostatic} (or \emph{isotropic}) part whereas the second term represents the \emph{deviatoric} part. The functions $k(\bx)$ and $\mu(\bx)$ denote the bulk and shear moduli, respectively, and $\theta(\bx)$ denotes the volumetric change and is given by 
\begin{equation}
 \theta(\bx) = \frac{3}{m} \Big( \int_{B_\varepsilon(\bx)} w(|\bz -\bx|)|\bz-\bx|\, |\bu(\bz) + \bz - (\bu(\bx) + \bx)| \,d\bz -m\Big),
\end{equation}
The radial function $w$ is given by 
\begin{equation}
 w(|\xi|) = \left\{\begin{aligned}
                    \frac{1}{|\bxi|^r} \qquad& \mbox{if $|\bxi|<\delta$} \\
                    0                 \qquad& \text{otherwise}
                   \end{aligned}
\right.
\end{equation}
and 
\begin{equation}
 m = \int_{\Omega} w(|\bxi|) |\bxi|^2 d\xi = \int_{B_{\delta}(0)} |\bxi|^{2-r} d\bxi = 4\pi\frac{\delta^{5-r}}{5-r} \qquad\text{for $r<5$}.
\end{equation}
Note that when $r<5$, $m$ is finite. For example, when $r=2$, $m = \frac{4}{3}\pi\delta^3 = |B_{\delta}(0)|$.

Let $\betab(\bx,\by)=\bu(\by) - \bu(\bx)$ denote the relative displacement. We linearize the peridynamic equation of motion with respect to small relative displacements, i.e., for $|\betab|\ll1$, as discussed  in \cite{Silling2010} . Observe that, in terms of $\betab$,
$$
\theta(\bx) = \frac{3}{m}\int_{B_\varepsilon(0)} w(|\bzeta|)|\bzeta|
(|\betab+\bzeta|-|\bzeta|)\,d\bzeta
$$
and thus
$$
  \frac{\partial\theta}{\partial\eta_i}=
  \frac{3}{m}\int_{B_\varepsilon(0)} w(|\bzeta|)|\bzeta|
\frac{\eta_i+\zeta_i}{|\betab+\bzeta|}\,d\bzeta
\qquad\mbox{for $i=1,2,3$}
$$
so that
$$
  \nabla_{\betab}\theta|_{\betab={\bf0}} = 
  \frac{3}{m}\int_{B_\varepsilon(0)} w(|\bzeta|)
\bzeta\,d\bzeta.
$$
Then, because $\theta=0$ when $\betab={\bf0}$, we have that
$$
  \theta_{lin}(\bx) =
  \frac{3}{m}\int_{B_\varepsilon(0)} w(|\bzeta|) \bzeta\cdot\betab\,d\bzeta =
  \frac{3}{m}\int_{B_\varepsilon(0)} w(|\bz-\bx|) (\bz-\bx)\cdot\big( \bu(\bz) - \bu(\bx)\big)\,d\bz,
$$
where $ \theta_{lin}$ denotes $\theta$ linearized about $\betab={\bf 0}$. Similarly, we find that
\begin{equation}\label{siglin}
\begin{aligned}
  \ttt_{lin}(\bx,\by) =
  \frac{3}{m}k(\bx)&w(|\bxi|)|\bxi| \frac{3}{m}\int_{B_\varepsilon(0)} w(|\bzeta|) \bzeta\cdot\betab\,d\bzeta
\\[1ex] &+ 
\frac{15}{m}\mu(\bx) w(|\bxi|) \Big( \frac{\bxi\cdot\betab}{|\bxi|}
- \frac{|\bxi|}{m} \int_{B_\varepsilon(0)} w(|\bzeta|) \bzeta\cdot\betab\,d\bzeta \Big)
\end{aligned}
\end{equation}
and
\begin{equation}\label{betlin}
\begin{aligned}
  \bgamma_{lin}(\bx,\by) &=
  \frac{\bxi}{|\bxi|} +   \Big(\frac{1}{|\bxi|} {\bf I} 
  - \frac{\bxi\otimes\bxi}{|\bxi|^3}\Big)\betab
  \\&=
  \frac{\by-\bx}{|\by-\bx|} +   \Big(\frac{1}{|\by-\bx|} {\bf I} 
  - \frac{(\by-\bx)\otimes(\by-\bx)}{|\by-\bx|^3}\Big)
 \big( \bu(\by)-\bu(\bx)\big).
\end{aligned}
\end{equation}
Therefore, from \eqref{pd2}, \eqref{siglin}, and \eqref{betlin}, and after ignoring higher-order terms in $\betab$, the linearized force density is given by
\begin{equation}\label{lpd1}
\begin{aligned}
 \TTT_{lin}(\bx,\by)
 =\bigg[&\frac{15}{m}\mu(\bx) w(|\by-\bx|)\frac{(\by-\bx)\otimes(\by-\bx)}{|\by-\bx|^2}\bigg]\big( \bu(\by)-\bu(\bx)\big)
\\& +\bigg[\frac{9}{m^2} w(|\by-\bx|)
 \Big(k(\bx)-\frac53\mu(\bx)\Big)
 \\&\qquad 
\times\Big( \int_{B_\varepsilon(\bx)} 
w(|\bz-\bx|) (\bz-\bx)\cdot\big( \bu(\bz)-\bu(\bx)\big)\,d\bz\Big)\bigg] (\by-\bx).
\end{aligned}
\end{equation}

Let
\begin{equation}\label{lpd2}
  (\CL \bu)(\bx)
  = \int_{B_\varepsilon(\bx)} \Big(\TTT_{lin}(\bx,\by)-\TTT_{lin}(\by,\bx)  \Big)\,d\by.
\end{equation}
Then, the linearized peridynamic equation of motion for a heterogeneous, isotropic solid is given by
$$
   \rho \ddot\bu = \CL \bu + \bb.
$$
The substitution of \eqref{lpd1} into \eqref{lpd2} yields
\begin{equation}\label{lpd3}
\begin{aligned}
(&\CL \bu)(\bx)
\\&  = \int_{B_\varepsilon(\bx)}
  \bigg[\frac{15}{m}\big(\mu(\bx)+\mu(\by)\big) w(|\by-\bx|)\frac{(\by-\bx)
\otimes(\by-\bx)}{|\by-\bx|^2}\bigg]\big( \bu(\by)-\bu(\bx)\big)
  \,d\by
  \\&\qquad
  + \int_{B_\varepsilon(\bx)}\int_{B_\varepsilon(\bx)}
  \bigg[\frac{9}{m^2}
 \Big(k(\bx)-\frac53\mu(\bx)\Big) w(|\by-\bx|)w(|\bz-\bx|)
 (\by-\bx)
\\&\qquad\qquad\qquad\qquad\qquad\qquad
\otimes(\bz-\bx)\big( \bu(\bz)-\bu(\bx)\big)\Big]
 \,d\bz d\by
\\&\qquad
  + \int_{B_\varepsilon(\bx)}\int_{B_\varepsilon(\by)}
  \bigg[\frac{9}{m^2}
 \Big(k(\by)-\frac53\mu(\by)\Big) w(|\by-\bx|)w(|\bz-\by|)
 (\by-\bx)
\\&\qquad\qquad\qquad\qquad\qquad\qquad
\otimes(\bz-\by)\big( \bu(\bz)-\bu(\by)\big)\Big]
 \,d\bz d\by.
\end{aligned}
\end{equation}

The following proposition shows that $(\CL \bu)(\bx)$ can be written as an integral operator acting on the relative displacement $\bu(\by)-\bu(\bx)$ for some kernel $\CCC$.

\begin{prop}\label{proplp}
The operator $\CL$ given by \eqref{lpd3} can be written as
\begin{equation}\label{lpd4}
(\CL \bu)(\bx)= \int_\Omega \CCC(\bx,\by)
\big( \bu(\by)-\bu(\bx)\big)\,d\by,
\end{equation}
where
$$
   \CCC(\bx,\by) = \KKK(\bx,\by) + \SSS(\bx,\by)
$$
with
\begin{equation}\label{lpd5}
   \KKK(\bx,\by) = \big(c_1(\bx)+c_1(\by)\big)
   w(|\by-\bx|)\frac{(\by-\bx)\otimes(\by-\bx)}{|\by-\bx|^2}
   \chi_{B_\varepsilon(\bx)}(\by)
\end{equation}
and
\begin{equation}\label{lpd6}
\begin{aligned}
&\SSS(\bx,\by) = \int_\Omega \Big[
c_2(\bz) w(|\bz-\bx|)w(|\by-\bz|)(\bz-\bx)\otimes(\by-\bz)
\chi_{B_\varepsilon(\bx)}(\bz)\chi_{B_\varepsilon(\by)}(\bz)
\\&\qquad
- 
c_2(\by) w(|\by-\bx|)w(|\bz-\by|)(\by-\bx)\otimes(\bz-\by)
\chi_{B_\varepsilon(\bx)}(\by)\chi_{B_\varepsilon(\by)}(\bz)
\\&\qquad
+ 
c_2(\bx) w(|\bz-\bx|)w(|\by-\bx|)(\bz-\bx)\otimes(\by-\bx)
\chi_{B_\varepsilon(\bx)}(\by)\chi_{B_\varepsilon(\bx)}(\bz)
\Big]\,d\bz,
\end{aligned}
\end{equation}
where
\begin{equation}\label{lpd7}
  c_1(\bx) = \frac{15}{m}\mu(\bx),\qquad
  c_2(\bx) = \frac{9}{m^2}  \Big(k(\bx) - \frac53\mu(\bx)\Big),
\end{equation}
and $\chi_{B_\varepsilon(\bx)}$ denotes the indicator function of the set ${B_\varepsilon(\bx)}$.
\end{prop}

{\em Proof.} It is obvious, with $c_1$ given by \eqref{lpd7}, that the first term in \eqref{lpd3} is equal to $\int_\Omega \KKK(\bx,\by)
\big( \bu(\by)-\bu(\bx)\big)\,d\by$. Thus, it remains to show that $\int_\Omega \SSS(\bx,\by)\big( \bu(\by)-\bu(\bx)\big)\,d\by$ is equal to the sum of the second and third terms in \eqref{lpd3}.

We first write
\begin{equation}\label{lpd8}
\int_\Omega \SSS(\bx,\by)\big( \bu(\by)-\bu(\bx)\big)\,d\by
= \int_\Omega \SSS(\bx,\by) \bu(\by)\,d\by
-
\int_\Omega \SSS(\bx,\by)\bu(\bx)\,d\by.
\end{equation}
For the first term in \eqref{lpd8}, we use \eqref{lpd6} to obtain
\begin{equation}\label{lpd9}
\begin{aligned}
&\int_\Omega \SSS(\bx,\by)\bu(\by)\,d\by 
\\&\quad
= \int_\Omega\int_\Omega \Big[
c_2(\bz) w(|\bz-\bx|)w(|\by-\bz|)(\bz-\bx)
\\&\qquad\qquad\qquad\qquad\qquad
\otimes(\by-\bz)
\chi_{B_\varepsilon(\bx)}(\bz)\chi_{B_\varepsilon(\by)}(\bz)\Big]\bu(\by)
\,d\bz d\by
\\&\qquad
-  \int_\Omega\int_\Omega \Big[
c_2(\by) w(|\by-\bx|)w(|\bz-\by|)(\by-\bx)
\\&\qquad\qquad\qquad\qquad\qquad
\otimes(\bz-\by)
\chi_{B_\varepsilon(\bx)}(\by)\chi_{B_\varepsilon(\by)}(\bz)\Big]\bu(\by)
\,d\bz d\by
\\&\qquad
+ \int_\Omega\int_\Omega  \Big[
c_2(\bx) w(|\bz-\bx|)w(|\by-\bx|)(\bz-\bx)
\\&\qquad\qquad\qquad\qquad\qquad
\otimes(\by-\bx)
\chi_{B_\varepsilon(\bx)}(\by)\chi_{B_\varepsilon(\bx)}(\bz)\Big]\bu(\by)
\,d\bz d\by
\\&\quad
= \int_{B_\varepsilon(\bx)}\int_{B_\varepsilon(\by)} \Big[
c_2(\by) w(|\by-\bx|)w(|\bz-\by|)(\by-\bx)\otimes(\bz-\by)
\Big]\bu(\bz)
\,d\bz d\by
\\&\qquad
-  \int_{B_\varepsilon(\bx)}\int_{B_\varepsilon(\by)}\Big[
c_2(\by) w(|\by-\bx|)w(|\bz-\by|)(\by-\bx)\otimes(\bz-\by)
\Big]\bu(\by)
\,d\bz d\by
\\&\qquad
+ \int_{B_\varepsilon(\bx)}\int_{B_\varepsilon(\bx)} \Big[
c_2(\bx) w(|\by-\bx|)w(|\bz-\bx|)(\by-\bx)\otimes(\bz-\bx)
\Big]\bu(\bz)
\,d\bz d\by,
\end{aligned}
\end{equation}
where, for the last equality, $\by$ and $\bz$ have been switched in the first and third integrals. Similarly, the second term in \eqref{lpd8}, after using \eqref{lpd6} and an appropriate change of variables, can be written as
\begin{equation}\label{lpd10}
\begin{aligned}
&\int_\Omega \SSS(\bx,\by)\bu(\bx)\,d\by 
\\&\quad
= \int_{B_\varepsilon(\bx)}\int_{B_\varepsilon(\by)} \Big[
c_2(\by) w(|\by-\bx|)w(|\bz-\by|)(\by-\bx)\otimes(\bz-\by)
\Big]\bu(\bx)
\,d\bz d\by
\\&\qquad
-  \int_{B_\varepsilon(\bx)}\int_{B_\varepsilon(\by)}\Big[
c_2(\by) w(|\by-\bx|)w(|\bz-\by|)(\by-\bx)\otimes(\bz-\by)
\Big]\bu(\bx)
\,d\bz d\by
\\&\qquad
+ \int_{B_\varepsilon(\bx)}\int_{B_\varepsilon(\bx)}\Big[ 
c_2(\bx) w(|\by-\bx|)w(|\bz-\bx|)(\by-\bx)\otimes(\bz-\bx)
\Big]\bu(\bx)
\,d\bz d\by
\\&\quad=
\int_{B_\varepsilon(\bx)}\int_{B_\varepsilon(\bx)}\Big[ 
c_2(\bx) w(|\by-\bx|)w(|\bz-\bx|)(\by-\bx)\otimes(\bz-\bx)
\Big]\bu(\bx)
\,d\bz d\by.
\end{aligned}
\end{equation}
The substitution of \eqref{lpd9} and \eqref{lpd10} into \eqref{lpd8} results in
$$
\begin{aligned}
&\int_\Omega \SSS(\bx,\by)\big( \bu(\by)-\bu(\bx)\big)\,d\by
\\
&= \int_{B_\varepsilon(\bx)}\int_{B_\varepsilon(\by)} \Big[
c_2(\by) w(|\by-\bx|)w(|\bz-\by|)(\by-\bx)\otimes(\bz-\by)
\Big]\big(\bu(\bz)-\bu(\by)\big)
\,d\bz d\by
\\
&+ \int_{B_\varepsilon(\bx)}\int_{B_\varepsilon(\bx)}\Big[ 
c_2(\bx) w(|\by-\bx|)w(|\bz-\bx|)(\by-\bx)\otimes(\bz-\bx)
\Big]\big(\bu(\bz)-\bu(\bx)\big)
\,d\bz d\by.\end{aligned}
$$
which, with $c_2$ given by \eqref{lpd7}, is equal to the sum of the last two terms in \eqref{lpd3}.\hfill$\Box$

\subsection{Relation between the peridynamics operator and the nonlocal operators}\label{pdandno}

Let
\begin{equation}\label{lpdw}
    w(|\bz|) = \frac{1}{|\bz|^2}
\end{equation}
and
$$
    \alpha(\bx,\by) = (\by-\bx)w(|\by-\bx|) \chi_{B_\varepsilon(\bx)}(\by) = \frac{\by-\bx}{|\by-\bx|^2}\chi_{B_\varepsilon(\bx)}(\by).
$$
Note that $\alpha(\bx,\by) = - \alpha(\by,\bx)$ and that the specialized Schwartz divergence kernel \eqref{specker} is given by
$$
   \brho_{\balpha}(\bx,\by,\bz) =
   \delta(\bx-\bz) \frac{\by-\bx}{|\by-\bx|^2}\chi_{B_\varepsilon(\bx)}(\by)
   +
   \delta(\bx-\by) \frac{\bz-\bx}{|\bz-\bx|^2}\chi_{B_\varepsilon(\bx)}(\bz).
$$

\begin{thm}\label{thmpdandno}
The linear peridynamic operator $\CL$ is given it terms of the operators of the nonlocal vector calculus by
\begin{equation}\label{lpnc1}
  -\CL u = 
  \CG_{\balpha}(c_1\CG^\ast_{\balpha} u)
  +\CG_{\balpha}(c_2\overline{\CG}^\ast_{\balpha} u)
\end{equation}
or, equivalently, by
\begin{equation}\label{lpnc2}
  -\CL u = 
  \CD_{\balpha}\big(c_1(\CD^\ast_{\balpha} u)^T\big)
  +\CG_{\balpha}(c_2\overline{\CG}^\ast_{\balpha} u).
\end{equation}
\end{thm}

{\em Proof.}
We observe that 
\begin{equation}\label{lpnc3}
\begin{aligned}
 \CG_{\balpha}(c_1\CG^\ast_{\balpha} \bu)
&= \int_{\bbR^3} \Big[c_1(\by)(\CG^\ast_{\balpha}\bu)(\by,\bx)
+c_1(\bx)(\CG^\ast_{\balpha}\bu)(\bx,\by)\Big]\balpha(\bx,\by)\,d\by
\\&=
-\int_{\bbR^3} \Big[c_1(\by)\big(\bu(\bx)-\bu(\by)\big)\cdot\balpha(\by,\bx) 
\\&
\qquad\qquad + c_1(\bx)\big(\bu(\by)-\bu(\bx)\big)\cdot\balpha(\bx,\by)\Big]\balpha(\bx,\by)\,d\by
\\&=
-\int_{\bbR^3} \big(c_1(\bx)+c_1(\by)\big)\Big[\big(\bu(\by)-\bu(\bx)\big)\cdot\balpha(\bx,\by)\Big] \balpha(\bx,\by)\,d\by
\\&=
-  \int_{\bbR^3} \big(c_1(\bx)+c_1(\by)\big)
\balpha(\bx,\by)\otimes \balpha(\bx,\by)
\big(\bu(\by)-\bu(\bx)\big)\,d\by
\\&=
-  \int_{B_\varepsilon(\bx)} \big(c_1(\bx)+c_1(\by)\big)
\frac{(\by-\bx)\otimes(\by-\bx)}{|\by-\bx|^4}
\big(\bu(\by)-\bu(\bx)\big)\,d\by.
\end{aligned}
\end{equation}
 
We next observe that
\begin{equation}\label{lpnc4}
\begin{aligned}
\CG_{\balpha}(c_2\overline{\CG^\ast_{\balpha}} \bu)
&= \int_{\bbR^3}
 \Big[c_2(\by)(\overline{\CG^\ast_{\balpha}}\bu)(\by)
+c_2(\bx)(\overline{\CG^\ast_{\balpha}}\bu)(\bx)\Big]\balpha(\bx,\by)\,d\by
\\&=
-\int_{\bbR^3}
 \bigg[c_2(\by)\int_{\bbR^3} \big(\bu(\bz)-\bu(\by)\big)
 \cdot \balpha(\by,\bz)\,d\bz
\\&
\qquad
 + c_2(\bx)\int_{\bbR^3} \big(\bu(\bz)-\bu(\bx)\big)
 \cdot \balpha(\bx,\bz)\,d\bz
 \bigg]\balpha(\bx,\by)\,d\by
\\&=
-\int_{\bbR^3}\int_{\bbR^3}
c_2(\by)\balpha(\bx,\by)\otimes \balpha(\by,\bz)\big(\bu(\bz)-\bu(\by)\big)\,d\bz d\by
\\&
\qquad
-\int_{\bbR^3}\int_{\bbR^3}
c_2(\bx)\balpha(\bx,\by)\otimes \balpha(\bx,\bz)\big(\bu(\bz)-\bu(\bx)\big)\,d\bz d\by
\\&=
-\int_{B_\varepsilon(\bx)} \int_{B_\varepsilon(\by)}
c_2(\by)\frac{(\by-\bx)}{|\by-\bx|^2}\otimes\frac{(\bz-\by)}{|\bz-\by|^2}
\big(\bu(\bz)-\bu(\by)\big)\,d\bz d\by
\\&
\qquad
-\int_{B_\varepsilon(\bx)}\int_{B_\varepsilon(\bx)}
c_2(\bx)\frac{(\by-\bx)}{|\by-\bx|^2}\otimes\frac{(\bz-\bx)}{|\bz-\bx|^2}
\big(\bu(\bz)-\bu(\bx)\big)\,d\bz d\by.
\end{aligned}
\end{equation}
Then, with $\CL\bu$  given by \eqref{lpd3}, $c_1(\bx)$ and $c_2(\bx)$ given by \eqref{lpd6}, and $w(\bx)$ given by \eqref{lpdw}, \eqref{lpnc1} follows from \eqref{lpnc3} and \eqref{lpnc4}.

Finally, \eqref{lpnc2} follows from the following proposition.
\hfill$\Box$

\begin{prop}
The operators $\CD$ and $\CG$ satisfy
$$
\CD_{\balpha}\big(c(\CD^\ast_{\balpha} \bu)^T\big)
=\CG_{\balpha}(c\CG^\ast_{\balpha} \bu)\qquad\mbox{for all $\bu$}.
$$
\end{prop}

{\em Proof.}
Using the definitions of $\CD_{\balpha}$, $\CD^\ast_{\balpha}$, $\CG_{\balpha}$, and $\CG^\ast_{\balpha}$, one finds
$$
\begin{aligned}
\CD_{\balpha}\big(c&(\CD^\ast_{\balpha} \bu)^T\big)
=
\int_{\bRn}
\Big(c(\by)(\CD^\ast_{\balpha}\bu)^T(\by,\bx)
+c(\bx)(\CD^\ast_{\balpha}\bu)^T(\bx,\by)\Big)
\balpha(\bx,\by)\, d\by
\\&=
-\int_{\bRn} \Big[ c(\by) \Big(\big(\bu(\bx)-\bu(\by)\big)\otimes\balpha(\by,\bx)\Big)^T
\\&\qquad\qquad
+c(\bx) \Big(\big(\bu(\by)-\bu(\by)\big)\otimes\balpha(\bx,\by)\Big)^T
\Big]\balpha(\bx,\by)\, d\by
\\&=
-\int_{\bRn} \big( c(\by) + c(\bx)\big)\Big(\balpha(\bx,\by)\otimes\big(\bu(\by)-\bu(\bx)\big) \Big)  \balpha(\bx,\by)\, d\by   
\\&=
-\int_{\bRn} \big( c(\by) + c(\bx)\big)  \Big(\big(\bu(\by)-\bu(\bx)\big)\cdot\balpha(\bx,\by)\Big) \balpha(\bx,\by)\, d\by
\\&=
-\int_{\bRn} \Big[ c(\by)\big(\bu(\bx)-\bu(\by)\big)\cdot\balpha(\by,\bx)
\\&\qquad\qquad
+ c(\bx)\big(\bu(\by)-\bu(\bx)\big)\cdot\balpha(\bx,\by) \Big]\balpha(\bx,\by)\,d\by
\\&=
\int_{\bRn} \Big[ c(\by)(\CG^\ast_{\balpha}\bu)(\by,\bx) +  c(\bx)(\CG^\ast_{\balpha}\bu)(\bx,\by)   \Big]\balpha(\bx,\by)\,d\by
= 
\CG_{\balpha}(c\CG^\ast_{\balpha} \bu).
\,\hfill\Box
\end{aligned}
$$

\section{Concluding remarks}
\label{sec_conclusion}

We close the paper by briefly considering a second intriguing special case of the nonlocal operators and also briefly discussing, in general terms, the types of situations in which the generalized operators may be of use.

\subsection{Another special case of the nonlocal operators}

A different simplification of the nonlocal operators is effected by the Schwartz kernel given by
\begin{equation}\label{speckeraa}
\bkappa_{\bbeta}(\bx,\by,\bz) = -\delta(\bx-\bz)\bbeta(\bx,\by) + \delta(\bx-\by)\bbeta(\bx,\bz)
\end{equation}
for a vector-valued function $\bbeta(\bx,\by)\in[L^1(\bRn\times\bRn)]^k$. It is easy to show that $\bkappa_{\bbeta}$ satisfies \eqref{divthm4} if and only if $\bbeta(\bx,\by)$ is {\em symmetric}, i.e., we have $\bbeta(\bx,\by) = \bbeta(\by,\bx)$. Note the contrasts with the simplified kernel $\bkappa_{\balpha}$ given by \eqref{specker} for which $\balpha$ is antisymmetric and the minus sign in \eqref{speckeraa} is replaced by a plus sign.

The specialized kernel $\bkappa_{\bbeta}$ is a divergence kernel, i.e., it satisfies \eqref{divthm4} and, for that kernel, the nonlocal operator $\CD$ is given by 
\begin{equation}\label{cdb}
 (\CD_{\bbeta}\bnu)(\bx) = \int_{\bRn}\big(\bnu(\by,\bx) - \bnu(\bx,\by) \big)\cdot\bbeta(\bx,\by)\,d\by.
\end{equation}
and the adjoint operator $\CD^\ast$ is given by
$$
 (\CD_{\bbeta}^\ast u)(\bx,\by) = -\big(u(\by) - u(\bx)\big)\bbeta(\bx,\by) .
$$
Note the difference in a sign between $\CD_{\balpha}$ and $\CD_{\bbeta}$ but the similarity in sign between $\CD_{\balpha}^\ast$ and $\CD_{\bbeta}\ast$. The other operators of the nonlocal calculus, e.g., $\CG$, $\CC$, and their adjoined operators, have the obvious definitions resulting from making or not making the sign changes in the definitions of the corresponding operators engendered by the kernel $\bkappa_{\balpha}$. In particular, the form of the nonlocal Laplacian operator remains unchanged, i.e., we have that
\begin{equation}\label{betalap}
  (\CD_{\bbeta}\CD_{\bbeta}^\ast u)(\bx) =-2
  \int_{\bRn} \big(u(\by) - u(\bx)\big)
  \bbeta(\bx,\by)\cdot \bbeta(\bx,\by)\,d\by.
\end{equation}

As was the case for the kernel $\bkappa_{\balpha}$, for the specialized Schwartz kernel $\bkappa_{\bbeta}$ we have that \eqref{ident1a}, \eqref{ident1b}, \eqref{ident1c}, \eqref{ident2c}, and \eqref{ident2e} hold. However, unlike the situation for $\bkappa_{\balpha}$, for $\bkappa_{\bbeta}$ we trivially have that \eqref{ident2a} also holds. On the other hand, the state-based peridynamic model of Section \ref{peridyn} cannot be expressed in terms of the nonlocal operators corresponding to the kernel $\bkappa_{\bbeta}$ as was done in Section \ref{pdandno} for the kernel $\bkappa_{\balpha}$. It would be of interest to explore what sort of mechanical model, if any, results from the use of the specialized operators $\CG_{\bbeta}$, $\CD_{\bbeta}$, etc., and to explore the differences between such a model and the state-based peridynamics model. Although beyond the scope of this work, this is a subject of current interest to the authors. However, an inkling of the differences can be gleaned by applying the operators $\CD_{\balpha}$ and $\CD_{\bbeta}$ to a scalar function. Conceptually, we do this by setting $\bnu(\bx,\by)={\bf a}u(\bx)$ in \eqref{divop} and \eqref{cdb}, where ${\bf a}$ denotes a constant vector, yielding
$$
(\CD_{\balpha}u)(\bx) = \int_{\bRn} \big(u(\by) + u(\bx)\big)\widehat\alpha(\bx,\by)\,d\by
$$
and
\begin{equation}\label{betalap2}
(\CD_{\bbeta}u)(\bx) =\int_{\bRn} \big(u(\by) - u(\bx)\big)\widehat\beta(\bx,\by)\,d\by,
\end{equation}
where $\widehat\alpha(\bx,\by)={\bf a}\cdot\balpha(\bx,\by)$ is an antisymmetric function and $\widehat\beta(\bx,\by)={\bf a}\cdot\bbeta(\bx,\by)$ is a symmetric function. Comparing the last equation in \eqref{mmmmmm} with \eqref{betalap2} after setting $\widehat\beta(\bx,\by)=-2\balpha(\bx,\by)\cdot\balpha(\bx,\by)$ and comparing \eqref{betalap} with \eqref{betalap2} after setting $\widehat\beta(\bx,\by)=-2\bbeta(\bx,\by)\cdot\bbeta(\bx,\by)$, we that the direct action of the operator $\CD_{\bbeta}$ on a scalar function $u(\bx)$ yields a nonlocal Laplacian of $u$.

\subsection{Modeling situations in which the general nonlocal operators can play a role}

Consider the specific Schwartz kernel given by
\begin{equation}\label{genkerab}
\bkappa_{{\lambda\balpha}}(\bx,\by,\bz) =  \lambda(\bx,\bz)\balpha(\bx,\by) + \lambda(\bx,\by)\balpha(\bx,\bz),
\end{equation}
where $\lambda(\bx,\by)$ and $\balpha(\bx,\by)$ are scalar and vector-valued functions, respectively. This kernel can be viewed  as a simplification of the general kernel $\bkappa$ or a generalization of the $\bkappa_{{\balpha}}$ given by \eqref{specker}. Clearly, if we set $\lambda(\bx,\by)=\delta(\bx-\by)$, then \eqref{genkerab} reduces to \eqref{specker}. Of course, we require the $\bkappa_{{\lambda\balpha}}(\bx,\by,\bz)$ to satisfy \eqref{divthm4} which implies that $\lambda(\bx,\by)$ and $\balpha(\bx,\by)$ must be such that
$$
    \int_{\bRn}\Big(\lambda(\bx,\bz)\balpha(\bx,\by) + \lambda(\bx,\by)\balpha(\bx,\bz)\Big) d\bx = {\bf 0}\qquad\forall\, \by,\bz\in\bRn.
$$

For the kernel \eqref{genkerab} we have, from \eqref{divop}, the nonlocal divergence operator $\CD_{{\lambda\balpha}}$ such that
\begin{equation}\label{genkerab1}
  (\CD_{{\lambda\balpha}}\bnu)(\bx) =
  \int_{\bRn}\int_{\bRn} 
  \Big(\lambda(\bx,\bz)\balpha(\bx,\by) + \lambda(\bx,\by)\balpha(\bx,\bz)\Big)\cdot\bnu(\by,\bz) \,d\bz d\by
\end{equation}
and, from \eqref{adjop}, the adjoint operator
\begin{equation}\label{genkerab2}
  (\CD_{{\lambda\balpha}}^\ast u)(\bx,\by) =
  \int_{\bRn} u(\bz) 
  \Big(\lambda(\bx,\bz)\balpha(\bx,\by) + \lambda(\bx,\by)\balpha(\bx,\bz)\Big)\, d\bz .
\end{equation}
Note that if $\lambda(\bx,\by)=\delta(\bx-\by)$, then \eqref{genkerab1} and \eqref{genkerab2} reduce to \eqref{sgrad} and\eqref{sgradadj}, respectively.

From \eqref{genkerab1}, we have that
$$
\begin{aligned}
(\CD_{{\lambda\balpha}}&\bnu)(\bx) 
\\&=
 \int_{\bRn}\int_{\bRn} 
  \lambda(\bx,\bz)\balpha(\bx,\by)\cdot\bnu(\by,\bz) \,d\bz d\by + 
  \int_{\bRn}\int_{\bRn}\lambda(\bx,\by)\balpha(\bx,\bz)\cdot\bnu(\by,\bz) \,d\bz d\by
  \\&=
  \int_{\bRn}\int_{\bRn} 
  \lambda(\bx,\bz)\balpha(\bx,\by)\cdot\bnu(\by,\bz) \,d\bz d\by + 
  \int_{\bRn}\int_{\bRn}\lambda(\bx,\bz)\balpha(\bx,\by)\cdot\bnu(\bz,\by) \,d\bz d\by
\\& = \int_{\bRn}\int_{\bRn} \big(\bnu(\bz,\by) + \bnu(\by,\bz)\big)\cdot\lambda(\bx,\bz)\balpha(\bx,\by) \,d\bz d\by
\end{aligned}
$$
so that
\begin{equation}\label{genkerab3}
(\CD_{{\lambda\balpha}}\bnu)(\bx) = 
\int_{\bRn} \balpha(\bx,\by)\cdot \bigg( \int_{\bRn} \big(\bnu(\bz,\by) + \bnu(\by,\bz)\big)\lambda(\bx,\bz) \,d\bz \bigg)d\by.
\end{equation}

\begin{figure}[h!]
\begin{center} 
\begin{tabular}{ccccc} 
\includegraphics[height=1.4in]{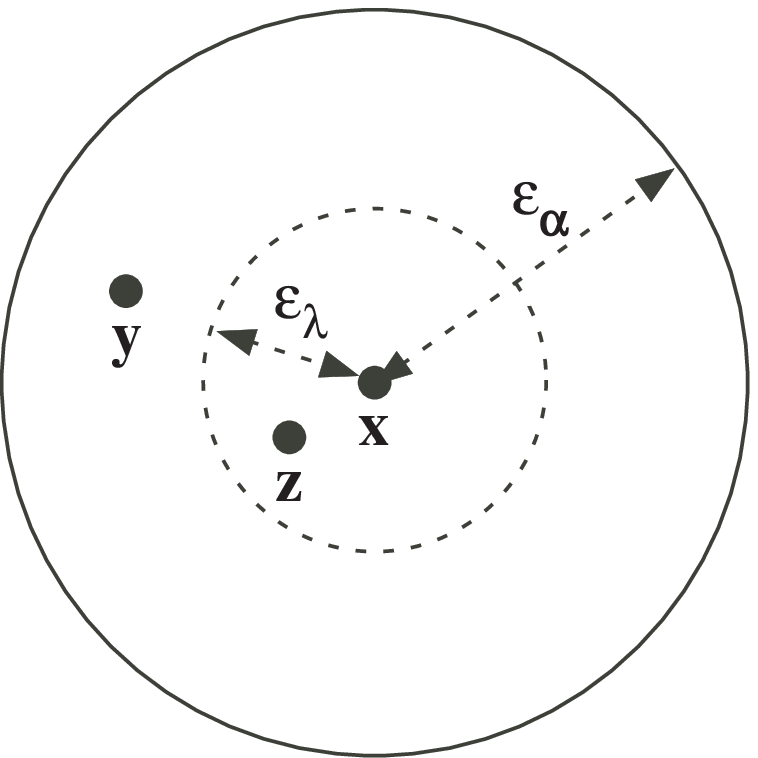}
&\raisebox{.7in}{\framebox{$\epsilon_\lambda\to0$}}
&\includegraphics[height=1.4in]{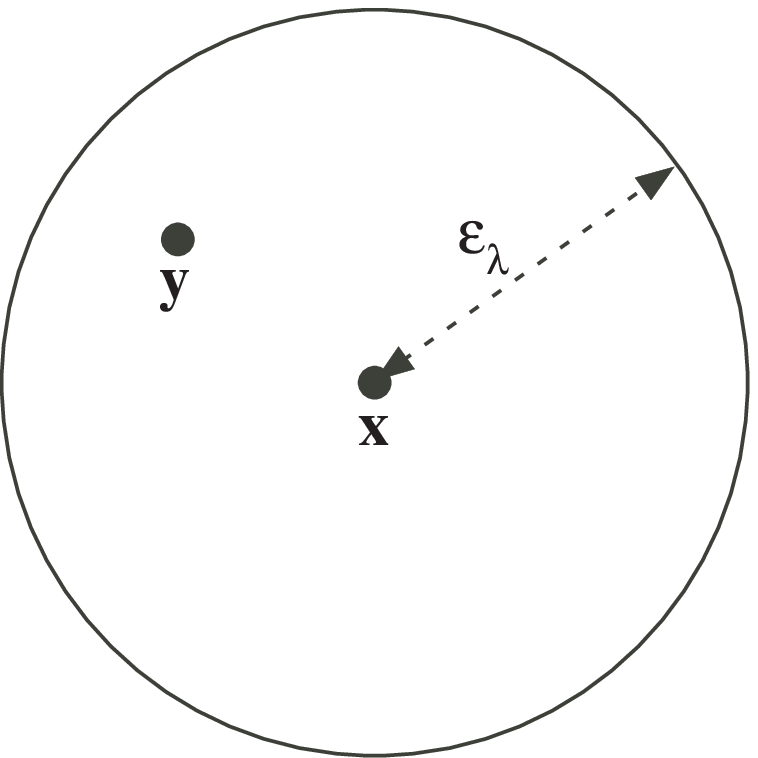}
&\raisebox{.7in}{\framebox{$\epsilon_{{\balpha}}\to0$}}
&\raisebox{.5in}{\includegraphics[height=.35in]{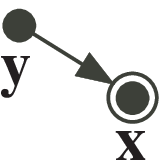}}
\end{tabular}
\caption{
Left: for the kernel \eqref{genkerab}, the contribution to the flux into a point $\bx$ from a point $\by$ in an $\epsilon_{{\balpha}}$-neighborhood of $\bx$ is determined by $\bz$ in an $\epsilon_\lambda$-neighborhood of $\bx$.
Middle: for the kernel \eqref{specker}, the contribution to the flux into a point $\bx$ from a point $\by$ in an $\epsilon_{{\balpha}}$-neighborhood of $\bx$ is determined only by the point $\by$.
Right: for local partial differential equation models, the contribution to the flux into a point $\bx$ is determined from points in an infinitesimal neighborhood of $\bx$.}
\label{fig1}
\end{center}
\end{figure}

Recall, from Section \ref{gennlc}, that the right-hand side of \eqref{genkerab3} is the total flux of $\bnu$ into the point $\bx$ coming from all points $\by\in\bRn$. Note that for the specialized kernel \eqref{specker}, i.e., for $\lambda(\bx,\by)=\delta(\bx-\by)$, the right-hand side \eqref{genkerab3} reduces to the right-hand side of \eqref{sgrad}. To better explain the difference between \eqref{genkerab3} and \eqref{sgrad}, it is instructive to consider the case of nonlocal interactions of finite extent, i.e., the case of $\lambda(\bx,\bz)$ and $\balpha(\bx,\by)$ having compact support. Specifically, we chose constants $\epsilon_\lambda$ and $\epsilon_{{\balpha}}$ such that $0< \epsilon_\lambda<\infty$ and $0<\epsilon_{{\balpha}} <\infty$ and then assume that
$$
\begin{aligned}
&\lambda(\bx,\bz) = 0  \qquad\mbox{for $\bz\not\in B_{\epsilon_\lambda}(\bx)$}
\\&
 \balpha(\bx,\by)={\bf0} \qquad\mbox{for $\by\not\in B_{\epsilon_{{\balpha}}}(\bx)$}.
\end{aligned}
$$
Then, the right-hand sides of \eqref{sgrad} and \eqref{genkerab3} become
\begin{equation}\label{genkerab4}
\mbox{\em flux into $\bx$} = \int_{B_{\epsilon_{{\balpha}}}(\bx)} \balpha(\bx,\by)\cdot\big(\bnu(\bx,\by) + \bnu(\by,\bx)\big) \,d\by
\end{equation}
and
\begin{equation}\label{genkerab5}
\mbox{\em flux into $\bx$} = \int_{B_{\epsilon_{{\balpha}}}(\bx)} \balpha(\bx,\by)\cdot \bigg( \int_{B_{\epsilon_\lambda}(\bx)} \big(\bnu(\bz,\by) + \bnu(\by,\bz)\big)\lambda(\bx,\bz) \,d\bz \bigg)d\by,
\end{equation}
respectively. Both \eqref{genkerab4} and \eqref{genkerab5} state that all points $\by$ in the ball centered at $\bx$ and of radius $\epsilon_{{\balpha}}$ contribute to the flux of $\bnu$ into the point $\bx$. However, \eqref{genkerab4} further states that the contribution to the flux of $\bnu$ into the point $\bx$ coming from a particular point $\by\in B_{\epsilon_{{\balpha}}}(\bx)$ only depends on the value\footnote{It is clear from \eqref{genkerab4} that only the symmetric part of $\bnu(\bx,\by)$ contributes to the flux so that there is actually no ambiguity between $\bnu(\bx,\by)$ and $\bnu(\by,\bx)$.} of $\bnu$ at the pair of points $\bx$ and $\by$. On the other hand, \eqref{genkerab5} states something quite different. We now have that the contribution to the flux of $\bnu$ into the point $\bx$ coming from a particular point $\by\in B_{\epsilon_{{\balpha}}}(\bx)$ depends on the values of $\bnu$ at all point pairs $\by$ and $\bz$ with $\bz\in B_{\epsilon_{{\balpha}}}(\bx)$. It is also clear that \eqref{genkerab5} reduces to \eqref{genkerab4} in an appropriate limit as $\epsilon_\lambda\to0$ and for an appropriate $\lambda(\bx,\bz)$, e.g., a Gaussian with variance and height depending $\epsilon_\lambda$ is such a way that it has unit area for all $\epsilon_\lambda$. We also know from previous work that, for appropriate $\alpha(\bx,\by)$, \eqref{genkerab4}, or more precisely $\CD_{{\balpha}}$ given by \eqref{sgrad}, reduces, as $\epsilon_{{\balpha}}\to0$, to the classical local differential divergence operator. We can interpret this limit as saying that only points in an infinitesimal ball centered at $\bx$ contribute to the flux into $\bx$, where an infinitesimal ball is needed so that one can be sure that derivatives are well defined. The sketches in Figure \ref{fig1} are meant to illustrate this discussion.

\end{document}